arXiv:0803.2119v1 [math.ST] 14 Mar 2008

# Jump estimation in inverse regression

**Leif Boysen and Axel Munk**

*L. Boysen, A. Munk*

*Institute for Mathematical Stochastics*

*Georg–August–Universität Göttingen*

*Maschmühlenweg 8–10*

*37073 Göttingen*

*Germany*

*E-mail:* boysen@math.uni-goettingen.de

munk@math.uni-goettingen.de

**Abstract:** We consider estimation of a step function $f$ from noisy observations of a deconvolution $\phi * f$, where $\phi$ is some bounded $L_1$-function. We use a penalized least squares estimator to reconstruct the signal $f$ from the observations, with penalty equal to the number of jumps of the reconstruction. Asymptotically, it is possible to correctly estimate the number of jumps with probability one. Given that the number of jumps is correctly estimated, we show that the corresponding parameter estimates of the jump locations and jump heights are $n^{-1/2}$ consistent and converge to a joint normal distribution with covariance structure depending on $\phi$, and that this rate is minimax for bounded continuous kernels $\phi$. As special case we obtain the asymptotic distribution of the least squares estimator in multiphase regression and generalisations thereof. In contrast to the results obtained for bounded $\phi$, we show that for kernels with a singularity of order $O(|x|^{-\alpha})$, $1/2 < \alpha < 1$, a jump location can be estimated at a rate of $n^{-1/(3-2\alpha)}$, which is again the minimax rate. We find that these rate do not depend on the spectral information of the operator rather on its localization properties in the time domain. Finally, it turns out that adaptive sampling does not improve the rate of convergence, in strict contrast to the case of direct regression.

## 1. Introduction

Assume we have observations from a regression model given by

$$Y = \Big( (\Phi f)\big(x_i\big) + \varepsilon_i \Big)_{i=1}^{n} , \tag{1}$$

where $\Phi f = \phi * f$ denotes convolution of some $L_1$-functions $\phi$ and $f$ and $\varepsilon_1, \varepsilon_2, \dots$ are i.i.d. mean zero random variables with finite second moment. In the following we denote model (1) as inverse (deconvolution) regression model and we assume throughout that $\phi$ is known. Suppose the objective function $f : [0,1] \to \mathbb{R}$ is in $L_1$ and moreover locally constant, i.e. a piecewise constant function with $k$ jumps given by

$$f(x) = \sum_{i=1}^{k+1} b_i 1_{[\tau_{i-1},\tau_i)}(x) , \tag{2}$$

s.t. $-\infty = \tau_0 \leq 0 < \tau_1 < \dots < \tau_k < 1 \leq \tau_{k+1} = \infty$ and $k \in \mathbb{N}$ possibly unknown (see Figure 1). From Figure 1 the difficulty of estimating jumps in inverse reression becomes visible: Due to the smoothing by $\phi$ jumps only appear as small changes in $\Phi f$.

In this paper we show that the joint least squares estimator $\hat{\theta}_n$ of jumps and heights

$$\theta = (b_1, \tau_1, b_2, \tau_2, \dots, b_k, \tau_k, b_{k+1}) \tag{3}$$

is $n^{-1/2}$ consistent and follows a multivariate normal limit law. This is in strict contrast to the case of direct regression (where $\Phi$ in (1) is the identity). In the latter case it is known that the LSE converges at the (minimax) $n^{-1}$ rate and its distribution (after recentering and rescaling with $n$) is given as the minimizer of a certain random walk process. Further, jump heights and locations are asymptotically independet (see van de Geer (1988); Yao and Au (1989); Müller (1992); Müller and Stadtmüller (1999); Yakir et al. (1999); Birgé and Massart (2006) for some references on jump estimation in direct regression). Finally, by an adaptive choice of the design points it is possible to speed up the $n^{-1}$ rate to any polynomial rate of convergence (Lan et al. (2007)). We will see that in inverse regression the situation is completely different w.r.t. all of these issues: In general, all components of $n^{1/2}(\hat{\theta}_n - \theta)$ will be dependent asymptotically (depending on the kernel $\phi$). Further, rather surprisingly, the $n^{-1/2}$ rate does not depend on the decay of the Fourier transform of $\phi$ which usually determines the rate of convergence in more common function spaces, such as Sobolev spaces (cf. Cavalier and Tsybakov (2002) among others). Indeed, we will show that the $n^{-1/2}$ rate is minimax if $\phi$ is a bounded, continuous function. Because our minimax lower bound will be independent of the design points we obtain the suprising finding that adaptive sampling cannot improve the rate of convergence in the inverse case.

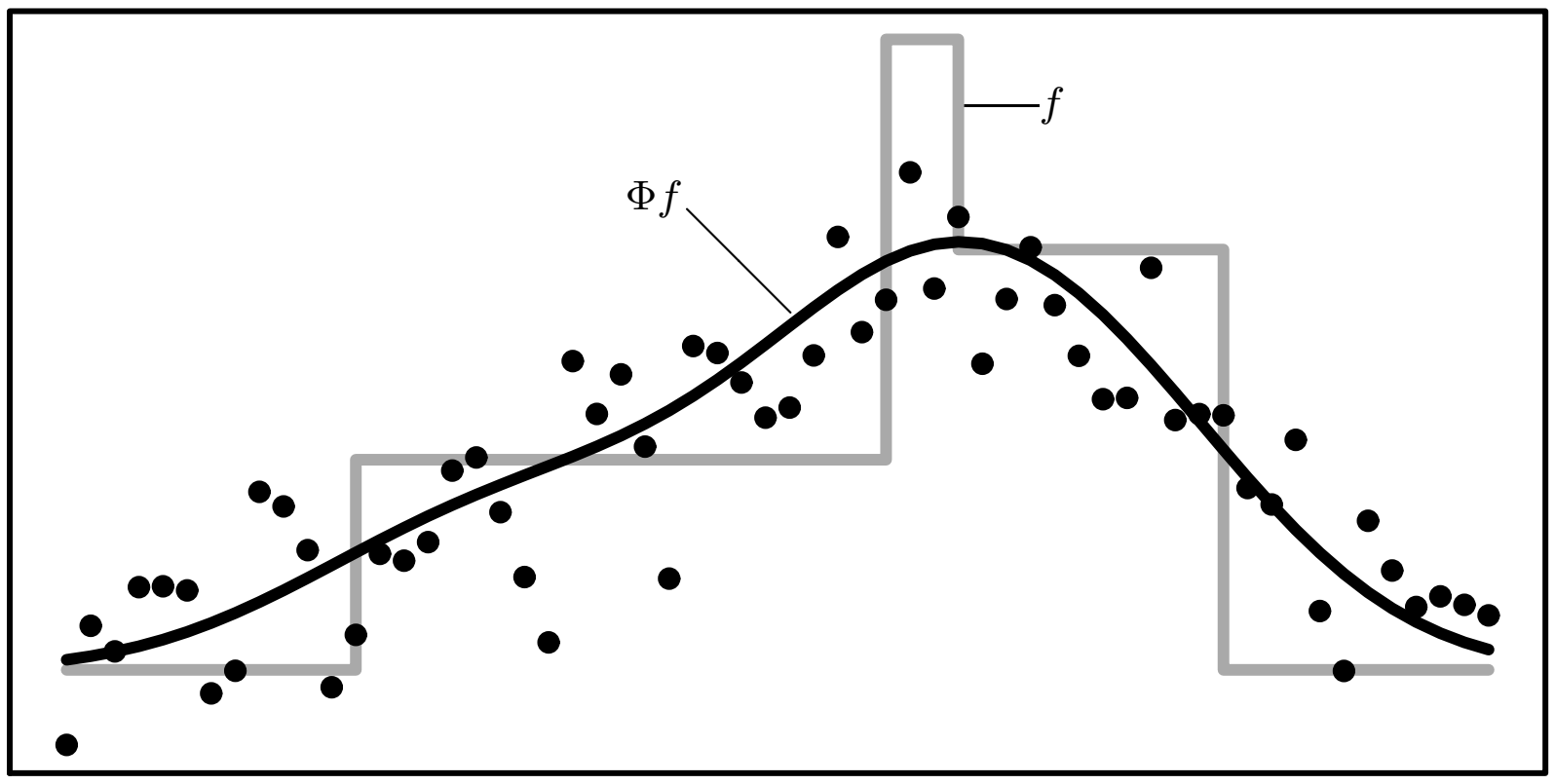

Fig 1. *Noisy observations of a blurred step function. The dots represent the observations and the black line the blurred function $\Phi f$, where $\Phi$ represents convolution with the gauss kernel. The gray line shows the original step function $f$, which is to be estimated.*

In fact a main motivation to consider the space of locally constant functions as in (2) stems from the observation that in general deconvolution is a difficult problem, which is reflected by rates of convergence which can be arbitrarily slow, e.g. $(\log n)^{-\beta}$ rates as for supersmooth (e.g. gaussian) deconvolution (cf.

Butucea and Tsybakov (2007)). However, we stress that in many practical situations, gaussian deconvolution is still applied, leading to satisfactory results (see e.g. Bissantz et al. (2007) for an example in astrophysics). At a first glance this seems to be contradictory. However, often a minimax result leads to rather pessimistic view, in particular in large function classes such as Sobolev spaces are. Often, more restrictive modeling is possible and necessary to obtain reasonably good rates of convergence. In fact the space of locally constant functions as considered in this paper (albeit of dimension $\infty$) yields a $n^{-1/2}$ rate of convergence *generically* which renders deconvolution in this setting as a practically feasable task. In fact, in this case the correct (and finite) number of jumps will be estimated asymptotically, and the problem reduces to a (nonsmooth) nonlinear regression problem.

We will give general conditions, which are sufficient to deduce the $n^{-1/2}$ rate. These conditions are borrowed from the theory of radial basis functions in native Hilbert spaces and from total positivity. They cover super-smooth functions such as the Gauss-kernel, polynomial kernels $\phi(x) = x^p\, 1_{[0,1)}(x)$ with $p = 0, 1, \dots$ and continuous symmetric functions $\phi$ which have a Fourier transform with an at most polynomial decay, satisfying $C(1 + |x|^{n_0})^{-1}$ for some $n_0 \in \mathbb{N}$, $C > 0$.

If the number of jumps is unknown, we show that – under the additional assumption of subgaussian tails of the error distribution – the number of jumps can be asymptotically estimated correctly with probability one.

We mention that our results can also be shown for more general Fredholm integral operators of the type $\Phi f = \int K(x, y) f(y) dy$ with continuous kernel $K : [0, 1] \times \mathbb{R} \to \mathbb{R}$ (see Boysen, 2006). For reasons of simplicity and ease of notation we do not treat this case here.

A classical model which fits into our framework was given by Quandt (1958). He introduced a linear regression model which obeys two separate regimes and where the change-point is not known. This model is called two-phase regression and inference in this setting was studied by Quandt (1960), Sprent (1961),

Hinkley (1969) and more recently by van de Geer (1988), Yakir et al. (1999) and Koul et al. (2003), among others. If the objective function $f$ is assumed to be continuous, two-phase regression can be modeled by an inverse regression model with a polynomial kernel with $p = 0$, i.e. $\phi(x) = 1_{[0,1)}(x)$. In this setting the $n^{-1/2}$ rate and the asymptotic distribution were derived by Hinkley (1969) and – for more general segmented regression models – by Feder (1975). From the perspective of a statistical inverse problem their results are quite natural to understand: multiphase regression corresponds to estimation of a jump function in a noisy Volterra equation where the location of jumps correspond to the kinks of the multiphase regression function.

Our results generalize the known results on the estimation of the intersection in two phase regression to the case where the objective function has an arbitrary number of phases and is piecewise polynomial of order $p+1$, with $p$ continuous derivatives and a $(p+1)$-th derivative, which is a step function. For piecewise linear regression ($p = 1$) in a deconvolution context this problem occurs in rheology where the relaxation time spectrum has to be estimated from measurements of the dynamic moduli of materials (cf. Roths et al., 2000). Other applications stem from biophysics, where the ion-channel activity of lipid membranes are measured by impedance spectroscopy and the jump locations indicate different opening states (cf. Schmitt et al., 2006; Römer et al., 2004). We obtain the somewhat suprising result that the rate of estimating the change-point does not depend on $p$, whereas in general nonparametric regression settings, the convergence rates for estimating a jump in the $p$-th derivative become slower as $p$ grows (see Raimondo, 1998).

The first one to investigate the change-point problem in the framework of a statistical inverse problem was Neumann (1997), who considered the estimation of a change-point in a density deconvolution model $Y = X + \xi$ with known error density $f_\xi$. He treated the case that the density of $X$ is bounded, has one jump at $\tau$ and is Lipschitz continuous elsewhere. In this setting $\tau$ can be estimated at a rate of $\min(n^{-1/(2\beta+1)}, n^{-1/(\beta+3/2)})$, provided the tails of the Fourier transform $\widehat{f_\xi}(x)$ decrease at a rate of $|x|^{-\beta}$. Moreover, he proved that these rates are optimal in a minimax sense. This result was extended by Goldenshluger et al. (2006c) (in a white noise model) to classes of functions $f$ which can be written as a sum of a step function and a function with smooth $m$-th derivative. They showed that in this case the minimax rates are of order $\min(n^{-1/(2\beta+1)}, n^{-(m+1)/(2\beta+2m+1)})$. If the smooth part of the function of interest belongs to a Paley-Wiener class, they show that a rate of $\min(n^{-1/2}, n^{-1/(2\beta+1)})$ can be obtained up to a logarithmic factor. Their recent work (Goldenshluger et al., 2006a,b) generalize these results to a unifying framework of sequence space models covering delay and amplitude estimation, estimation of change-points in derivatives and change point estimation in a convolution white noise model. We remark that the specific choice of jump functions in (2) used in this work comes close to the super-smooth case for $\beta \geq 1/2$, but we can get rid of the additional logarithmic factor. Moreover, we will see that similar rates hold in the case of $\beta < 1/2$ if the assumption on the boundedness of the kernel is dropped (see Remark 3).

This work is structured as follows. Section 2 gives some basic notation and the main assumptions. The estimate and its asymptotic properties are given in section 3 and the proof of the main result can be found in section 4. In section 5 we derive the required results from the theory of radial basis functions whicdh yields sufficient conditions on $\phi$ for the asymptotic normality of the LSE. Finally, in section 6 we derive the minimax rate for estimating the jump location.

## 2. Model assumptions and Notation

### 2.1. Notation

Define

$$\Gamma_k := \{(\gamma_0, \gamma_1, \ldots, \gamma_{k+1}) \,:\, -\infty = \gamma_0 < 0 < \gamma_1 < \ldots < \gamma_k < 1 < \gamma_{k+1} = \infty\}$$

as the set of possible jumps of $f$ in (1), and denote the corresponding function space of locally constant functions with at most $k$ jumps by

$$T_k := \Big\{ \sum_{i=1}^{k+1} b_i 1_{[\tau_{i-1},\tau_i)}(x) \ : \tau \in \Gamma_k \, , \, b_i \in \mathbb{R} \Big\} \, .$$

Write $T_\infty := \bigcup_{k=1}^{\infty} T_k$ for the set of all step functions on $\mathbb{R}$ with a finite but arbitrary number of jumps, where we exclude an isolated jump at the end points of the interval $[0,1]$. Note, that outside of $[0,1]$ these functions are constant. Let $T_{k,R} = \{g \in T_k : \|g\|_\infty < R\}$ as well as $T_{\infty,R} := \bigcup_{k=1}^{\infty} T_{k,R}$ the corresponding spaces of uniformly bounded functions for some $R > 0$. If not mentioned otherwise, the restriction of these spaces to $[0,1]$ are considered to be subspaces of $L_2([0,1])$.

Define the empirical norm $\|\cdot\|_n$ and the empirical inner product $\langle \cdot, \cdot \rangle_n$ by

$$\|f\|_n^2 := \frac{1}{n} \sum_{i=1}^{n} f^2(x_i) \qquad \text{as well as} \qquad \langle f, g \rangle_n := \frac{1}{n} \sum_{i=1}^{n} f(x_i) g(x_i),$$

where $x_1, \ldots, x_n$ are the design points. Similarly set

$$\|y\|_n^2 := \frac{1}{n} \sum_{i=1}^{n} y_i^2 \qquad \text{as well as} \qquad \langle y, z \rangle_n := \frac{1}{n} \sum_{i=1}^{n} y_i z_i$$

for $y, z \in \mathbb{R}^n$.

Write $g(t_+) := \lim_{x \searrow t} g(x)$ for the right limit of $g$ in $t$ and $g(t_-) := \lim_{x \nearrow t} g(x)$ for the corresponding left limit. For some proper function $g : \mathbb{R} \to \mathbb{R}$ define the set of jump points of g as

$$\mathcal{J}(g) := \{t \in [0,1] : g(t_-) \neq g(t_+)\} \tag{4}$$

and $J_\#(f) := \#\mathcal{J}(f) + 1$, where $\#\mathcal{J}(f)$ denotes the number of jumps, which may be infinite.

Define the distance of some point $a \in \mathbb{R}$ to the set $B \subset \mathbb{R}$ as

$$\mathrm{d}(a, B) = \inf_{b \in B} |a - b|$$

and, slightly abusing notation, the Hausdorff distance of two sets $A, B$ as

$$\mathrm{d}(A, B) = \max\{\sup_{a \in A} \mathrm{d}(a, B) \, , \sup_{b \in B} \mathrm{d}(b, A)\} \, .$$

Finally, for ease of notation for any $a, b \in \mathbb{R}$, $[a, b]$ and $(a, b)$ always denote the intervals $[\min(a, b), \max(a, b)]$ and $(\min(a, b), \max(a, b))$, respectively.

### *2.2. Assumptions*

**Assumptions on the error** If the number of jumps is known the following basic assumption is sufficient to deduce the $n^{-1/2}$ rates of convergence for the least squares estimates.

**Assumption A.** *The array $(\varepsilon_1, \dots, \varepsilon_n)$ consists of independent identically distributed random variables with mean zero for every $n$. Additionally, assume*

$$\mathrm{E}(\varepsilon_1^2) = \sigma^2 < \infty \,.$$

If the number of jumps of the objective function is unknown, we will additionally need that the error satisfies the following subgaussian condition.

(A1) There exists some $\alpha > 0$ such that $\mathrm{E}(\exp(\varepsilon_1^2/\alpha)) < \infty$.

**Assumptions on the kernel** The parameters of $f$ in (1) and (2) are identifiable if

$$f \in T_{k,r} \quad \text{and} \quad 0 = \left\|(\Phi f)(\cdot)\right\|_2 \quad \Rightarrow \quad f \equiv 0. \tag{5}$$

By

$$\left\|(\Phi f)(\cdot)\right\|_2 = \left\| \sum_{i=1}^{k+1} b_i (\Phi 1_{[\tau_{i-1}, \tau_i)})(\cdot) \right\|_2$$

relation (5) holds if the functions

$$(\Phi 1_{[\tau_0, \tau_1)})(\cdot), (\Phi 1_{[\tau_1, \tau_2)})(\cdot), \dots, (\Phi 1_{[\tau_k, \tau_{k+1})})(\cdot) \tag{6}$$

are linearly independent in $L_2([0, 1])$ for any $(\tau_0, \dots, \tau_{k+1}) \in \Gamma_k$.

Throughout the following we require a slightly stronger condition, the independence of the functions in (6) together with their derivatives.

**Assumption B.** *Let*

$$\Delta_\phi(x,a,b) := \begin{cases} \int_a^b \phi(x-y)dy & b \neq a\,, \\ \\ \phi(x-a) & b = a\,. \end{cases} \tag{7}$$

*Assume that $\phi \in L_1(\mathbb{R}) \cap L_2(\mathbb{R}) \cap L_\infty(\mathbb{R})$ is piecewise continuous with finitely many jumps. Additionally the functions*

$$\Delta_\phi(x,\tau_0,\tau_1)\,,\ \Delta_\phi(x,\tau_1,\tau_2)\,,\ \ldots\,,\ \Delta_\phi(x,\tau_k,\tau_{k+1})$$

*are linearly independent for every choice of $k \in \mathbb{N}$ and*

$$-\infty = \tau_0 \leq 0 < \tau_1 \leq \tau_2 \leq \ldots \leq \tau_k < 1 \leq \tau_{k+1} = \infty\,,$$

*where only two subsequent $\tau_i$ are allowed to be equal.*

The following theorem gives some general conditions, which are sufficient for $\phi$ to satisfy Assumption B.

**Theorem 2.1.** *The function $\phi$ satisfies Assumption B if one of the following conditions is satisfied.*

*(i)* *$\phi \in C(\mathbb{R}) \cap L_1(\mathbb{R})$ is a symmetric real-valued function with Fourier transform $\widehat{\phi}(x) \geq 0$, such that there exists $n_0 \in \mathbb{N}$ and $C > 0$ with*

$$C(1+|x|^{n_0})^{-1} \leq |\widehat{\phi}(x)| \qquad \textit{for all}\, x \in \mathbb{R}\,. \tag{8}$$

*(ii)* *$\phi$ is extended sign regular of order $k+2$ on $\mathbb{R}$, with $0 < \int \phi(x)dx < \infty$.*

*(iii)* *The function $\phi$ is given by*

$$\phi(x) = \begin{cases} x^p & x \in [0,1] \\ \\ 0 & else \end{cases} \qquad p \in \{0,1,2,\ldots\}\,.$$

The proof of part (i) is given in section 5, the proofs of part (ii) and (iii) are straightforward and can be found in Boysen (2006), Section 5.2 and 5.3. A definition of sign-regularity can, for example, be found in Karlin and Studden (1966). Note that part (ii) covers the Gauss kernel $\phi(x) = (2\pi)^{-1/2}\exp(-x^2/2)$ (see Section 3, Example 5 in Karlin and Studden, 1966).

Examples of kernels which satisfy condition (i) are the Laplace kernel $\phi(x) = \exp(-|x|)/2$, the kernel $\phi(x) = \cos(x)\exp(-|x|)$ and kernels of the type $\phi(x) = (1-|x|)_+^p$ for $p = 2, 3, \ldots$ where $x_+$ denotes the positive part of $x$. Moreover, the convolution of any two kernels $\phi_1, \phi_2$ satisfying (i)-(iii) clearly also satisfies this condition.

**Assumptions on the design points** We make the following assumption on the design points.

**Assumption C.** *There exists a function* $h : [0,1] \to [c_l, c_u]$ *with* $0 < c_l < c_u < \infty$ *and* $\int_0^1 h(x)dx = 1$*, such that*

$$\frac{i}{n} = \int_0^{x_{(i)}} h(x)dx + \delta_i$$

*for all* $i = 1, \ldots, n$*, with* $\max_{i=1,\ldots,n} |\delta_i| = O_P(n^{-1/2})$*.*

*Moreover, the design points* $x_1, \ldots, x_n$ *are independent of the error terms* $\varepsilon_1, \ldots, \varepsilon_n$*. Here* $x_{(i)}$ *denotes the i-th order statistic of* $x_1, \ldots, x_n$*.*

Note that the above assumption covers random designs as well as fixed designs generated by a regular density in the sense of Sacks and Ylvisaker (1970). If the design points $x_1, \ldots, x_n$ are nonrandom, the $O_P(n^{-1/2})$ term above is to be understood as $O(n^{-1/2})$. In this case the design points have to be understood as a triangular scheme. Dümbgen and Johns (2004) use a similar assumption on the design points.

Note that if we assume a triangular scheme and fixed design points instead, all results can be obtained essentially in the same way. The only argument which has to be slightly modified, is the one based on the law of the iterated logarithm in the proof of Lemma 4.13. Note that, the respective inequalities remain valid because error terms in a triangular scheme can be replaced distributionally equivalent by a sequence of i.i.d. random variables.

## 3. Estimate and asymptotic results

**Estimate** Define the restricted least squares estimate $\hat{f}_n$ as approximate minimizer of the empirical $L_2$ distance to the data in the space $T_{k,R}$. More precisely $\hat{f}_n \in T_{k,R}$ and

$$\|\Phi \hat{f}_n - Y\|_n \leq \min_{g \in T_{k,R}} \left( \|\Phi g - Y\|_n^2 \right) + o_p(n^{-1})\,. \tag{9}$$

The minimizer of the functional on the right hand side always exists (compare Lemma 4.6). Note that we do not assume that the minimum is attained, but only that the functional above can be minimized up to some term of order $o_p(n^{-1})$. It does not need to be unique. This assumption allows for numerical approximation of the minimizer and gives an intuition of the needed precision for the asymptotic results to be valid. The restriction to functions with $\|f\|_\infty < R$ is a technical assumption, which requires that some upper bound of the supremum norm of the objective function is known beforehand.

Note that any estimator $\hat{f}_n$ has a representation as

$$\hat{f}_n(x) = \sum_{i=1}^{k+1} \hat{b}_i 1_{[\hat{\tau}_{i-1}, \hat{\tau}_i)}(x)\,, \tag{10}$$

with vectors $\hat{b} = (\hat{b}_1, \ldots, \hat{b}_{k+1})^t$ and $\hat{\tau} = (\hat{\tau}_0, \ldots, \hat{\tau}_{k+1})^t$, which are the approximate least squares estimates (in the sense of (9)) of the true parameter vectors $b$ and $\tau$ given by equation (2).

If the number of jumps is unknown, a different estimate is needed. In this case, assume that the penalized least squares estimate $\hat{f}_{\lambda_n}$ satisfies $\hat{f}_{\lambda_n} \in T_{\infty,R}$ and is defined as any solution of

$$\|\Phi \hat{f}_{\lambda_n} - Y\|_n + \lambda_n J_\#(\hat{f}_{\lambda_n}) \leq \min_{g \in T_{\infty,R}} \left( \|\Phi g - Y\|_n^2 + \lambda_n J_\#(g) \right) + o_p(n^{-1})\,, \tag{11}$$

where $\lambda_n > 0$ is some smoothing parameter, s.t. $\lambda_n \to 0$ as $n \to \infty$.

**Asymptotic results** Before we state the main result, we first define the map $\nu : [0,1] \mapsto \mathbb{R}^{2k+1}$ by

$$\nu(x) = \begin{pmatrix} \Delta_\phi(x, \tau_0, \tau_1) \\ (b_1 - b_2)\Delta_\phi(x, \tau_1, \tau_1) \\ \Delta_\phi(x, \tau_1, \tau_2) \\ \vdots \\ (b_k - b_{k+1})\Delta_\phi(x, \tau_k, \tau_k) \\ \Delta_\phi(x, \tau_k, \tau_{k+1})) \end{pmatrix}, \tag{12}$$

and the $(2k+1) \times (2k+1)$ matrix $V$ by its entries

$$(V)_{ij} = \int_0^1 (\nu(x)\nu(x)^t)_{ij}\, h(x)dx\,. \tag{13}$$

Here $h$ is the design density given by Assumption C. Now we are able to formulate the asymptotic result for the least squares estimator.

**Theorem 3.1.** *Suppose the Assumptions A, B and C are met. Let $\hat{f}_n$ and $V$ be given by (10) and (13), respectively. Set $\theta$ as the parameter vector of $f$ given in (3), and $\hat{\theta}_n$ as the corresponding vector of estimates defined by (10). Given (9) and model (1), then*

*(i)* $\sqrt{n}(\hat{\theta}_n - \theta) \xrightarrow{\mathcal{L}} N(0, \sigma^2 V^{-1})$.

*Moreover,*

*(ii)* $\|\Phi f - \Phi \hat{f}_n\|_2 = O_P(n^{-1/2})$.

*(iii)* $\mathrm{d}(\mathcal{J}(f), \mathcal{J}(\hat{f}_n)) = O_P(n^{-1/2})$.

*(iv)* $\|f - \hat{f}_n\|_2 = O_P(n^{-1/4})$.

*(v) $V$ is positive definite.*

The following theorem implies that the penalized and the restricted least squares estimates asymptotically coincide, i.e. the number of jumps in $T_\infty$ is asymptotically correctly estimated with probability one. In this sense the the results of Theorem 3.1 can be applied to the penalized estimate $\hat{f}_{\lambda_n}$.

**Theorem 3.2.** *Suppose condition (A1), (11) and the assumptions of Theorem 3.1 are satisfied. If* $\lambda_n \to 0$ *and* $\lambda_n n^{1/(1+\epsilon)} \to \infty$ *for some* $\epsilon > 0$ *as* $n \to \infty$, *then*

$$\lim_{n\to\infty} P\big(\#\mathcal{J}(\hat{f}_{\lambda_n}) = \#\mathcal{J}(f)\big) = 1\,.$$

The proofs of Theorem 3.1 and 3.2 can be outlined as follows. For a known number of jumps an entropy argument yields consistency of the least squares estimator. It is possible to represent the estimator as the minimizer of a stochastic process, which allows for a local stochastic expansion. This can be used to derive asymptotic normality. If the number of jumps is unknown, an imitation of techniques from empirical process theory shows that for a suitable choice of the smoothing parameter the case of an unknown number of jumps can asymptotically be reduced to the case where this number is known.

The details of the proofs are given in several steps in section 4.

The next theorem states that the rate given above is optimal in a minimax sense.

**Theorem 3.3.** *Suppose the Assumption B is met and* $\varepsilon_1, \ldots, \varepsilon_n$ *are independent identically distributed normal random variables with zero mean and positive variance. Set*

$$\Theta = \left\{\theta = (b_1, \tau_1, b_2, \tau_2, \ldots, b_k, \tau_k, b_{k+1}) : f_\theta(\cdot) := \sum_{i=1}^{k+1} b_i 1_{[\tau_{i-1},\tau_i)}(\cdot) \in T_{k,R}\right\}.$$

*For arbitrary fixed design points* $x_1, \ldots, x_n \in [0,1]$ *denote by* $P_\theta^n$ *the probability measure associated with the observations*

$$Y_i = (\Phi f_\theta)(x_i) + \varepsilon_i \quad i = 1, \ldots, n.$$

*Then there exists some* $c_0 > 0$ *independent of* $n$ *and* $x_1, \ldots, x_n$ *such that*

$$\inf_{\hat\theta} \sup_{\theta\in\Theta} \mathrm{P}_\theta^n(\|\hat\theta - \theta\|_2 \geq c_0 n^{-1/2}) > 0\,.$$

The proof is given in section 6.

**Remarks and Extensions**

**Remark 1.** *(Adaptive sampling). Theorem 3.3 states for any fixed bounded kernel $\phi$ and any choice of design points, that faster rates of convergence as $n^{-1/2}$ are not possible. This is intuitively clear as the convolution "spreads" the information of the jump location over the whole interval. As a consequence adaptive sampling schemes (where the sampling point $x_i$ may dependent on the data $Y_1, \cdots, Y_{i-1}$) cannot lead to a faster rate of convergence as $n^{-1/2}$. This is in strict contrast to the case of direct regression ($\Phi = Id$) where any polynomial rate of convergence can be achieved by an adaptice scheme (Lan et al. (2007)).*

**Remark 2.** *(Noisy Fredholm equations). All results of this chapter can also be shown for more general integral operators of the type $\Phi f = \int K(x,y)f(y)dy$ with continuous kernel $K : [0,1] \times \mathbb{R} \to \mathbb{R}$ satisfying $\sup_{x \in 0,1} \|K(x,\cdot)\|_{L_1} < \infty$. In this case in definition (7) $\phi(x-y)$ has to be replaced by $K(x,y)$. Assumption B can be formulated in the same way.*

**Remark 3.** *(Singular kernels). If the assumption of the boundedness of the integral kernel is dropped, faster rates than $O_P(n^{-1/2})$ for estimating the jump location can be achieved. Indeed if $\phi$ is an Abel type kernel $\phi_\alpha(x) = x^{-\alpha} 1_{(0,\infty)}(x)$ for $\alpha \in (0,1)$ then a jump can be recovered at a rate of $O_P(n^{-1/\min(2,3-2\alpha)})$. Given a uniform design, these rates are minimax. For details see Boysen (2006), chapter 8.2. We mention that in this case adaptive sampling can improve the rate of convergence similar as in the direct case but in contrast to a bounded kernel (cf. again Remark 1). Interestingly, the $n^{-1}$ sampling rate is achieved as $\alpha \to 1$, which is well known to be the best possible rate in direct regression for the estimation of a jump. Hence, singular kernels, with a spike at least as strong as $|x|^{-1}$ already localize jumps with the same rate as for the direct case, which is achieved as $\alpha \to \infty$.*

*Note that the "elbow" in the rates of convergence occurs at $\alpha = 1/2$, and that the $n^{-1/2}$ rate holds for the case where $\phi_\alpha$ is square integrable on bounded*

*intervals.*

*This corresponds to findings of Neumann (1997) and Goldenshluger et al. (2006c), who also observe an elbow in the rate of convergence of recovering a change point in an inverse problem at* $\beta = 1/2$*, if the Fourier transform of* $\widehat{\phi}(x)$ *decreases at rate of* $|x|^{-\beta}$*. Goldenshluger et al. (2006c) give a rate of* $O_P(n^{-1/\min(2,2\beta+1)})$ *up to a logarithmic term if the smooth part of the function of interest is in a Paley-Wiener class. From*

$$\left|\widehat{\phi}_\alpha(x)\right| = |x|^{-1+\alpha}\Gamma(1-\alpha)$$

*it follows that the "elbow" for* $\beta = 1/2$ *can be identified with the elbow for* $\alpha = 1/2$*.*

## 4. Proof of Theorem 3.1 and 3.2

We start with some technical lemmata, give some entropy results on the spaces of interest which are required to apply tools of empirical process theory to prove consistency of the estimates. Afterwards we give a local stochastic expansion of the minimized process and use this to derive asymptotic normality. Finally we again imitate some techniques from empirical process theory to show that the penalized estimate asymptotically coincides with the restricted least squares estimate. Note that Assumption B is needed to assure identifiability as well as positive definiteness of the asymptotic covariance matrix $V$.

### *4.1. Some technical lemmata*

In order to gain some insight into the model, it is useful to have a closer look at the implications of Assumption B on the mapping $\Phi$ restricted to the space of step functions. The following lemma collects some properties of this mapping.

**Lemma 4.1.** *Given Assumption B the following holds true.*

*(i) For all* $\epsilon > 0$ *there exists* $0 < C_0 < \infty$ *such that for all* $f \in T_\infty$

$$\|\Phi f\|_n^2 \le C_0 \|f\|_{L_2([-\epsilon,1+\epsilon])}^2 .$$

*(ii) For all* $\epsilon > 0$ *the map* $\Phi : (\, T_k, \|\cdot\|_{L_2([-\epsilon,1+\epsilon])}\,) \to L_2([0,1])$ *is continuous.*

*(iii)* $\Phi : T_k \to L_2([0,1])$ *is one-to-one.*

*(iv) The function* $(\Phi f)$ *is Lipschitz continuous on* $\mathbb{R}$ *for all* $f \in T_\infty$*.*

*Proof.* By Assumption B we have that $\|\phi\|_\infty = C < \infty$. Hence

$$\begin{aligned}\|\Phi f\|_n^2 &\le \int f^2(y) \frac{1}{n} \sum_{i=1}^n \phi^2(x_i - y) dy \\ &\le C^2 \int_0^1 f^2(y) dy + \int_{\mathbb{R}\setminus[0,1]} f^2(y) \frac{1}{n} \sum_{i=1}^n \phi^2(x_i - y) dy ,\end{aligned}$$

for $f \in T_\infty$. Note that $f$ is constant on $(-\infty, 0)$ and $[1, \infty)$. This gives

$$\begin{aligned}\|\Phi f\|_n^2 &\le C^2 \|f\|_2^2 + \|\phi\|_{L_2(\mathbb{R})}^2 (f|_{(-\infty,0)})^2 + \|\phi\|_{L_2(\mathbb{R})}^2 (f|_{[1,\infty)})^2 \\ &= C^2 \|f\|_2^2 + \frac{\|\phi\|_{L_2(\mathbb{R})}^2}{\epsilon} \Big( \int_{-\epsilon}^0 f^2(y) dy + \int_1^{1+\epsilon} f^2(y) dy \Big) \\ &\le C_0 \|f\|_{L_2([-\epsilon,1+\epsilon])}^2 ,\end{aligned}$$

for some $C_0$ depending on $\phi$ and $\epsilon$ only. This proves (i).

Similarly we can show $\|\Phi f\|_2 \le C \|f\|_{L_2([-\epsilon,1+\epsilon])}^2$ for $f \in T_k$ which gives continuity and hence (ii). As argued in the part on the assumptions on the kernel in section 2.2, (iii) follows from the independence of $\Delta_\phi(\cdot, \tau_i, \tau_{i+1})$.

To prove (iv), note that

$$\begin{aligned}|(\Phi 1_{[a,b)})(x) - (\Phi 1_{[a,b)})(x+\delta)| &= \Big| \int_{x-b}^{x-a} \phi(y) dy - \int_{x+\delta-b}^{x+\delta-a} \phi(y) dy \Big| \\ &\le 2|\delta| \|\phi\|_\infty ,\end{aligned}$$

for any $x, \delta \in \mathbb{R}$ and $a, b \in \mathbb{R} \cup \{-\infty, \infty\}$. For $f \in T_\infty$ with $\#\mathcal{J}(f) < \infty$, this gives $|(\Phi f)(x) - (\Phi f)(x+\delta)| \le |\delta| (2\#\mathcal{J}(f) \|f\|_\infty \|\phi\|_\infty)$. □

The following lemma provides a link of the empirical and the $L_2$ norm.

**Lemma 4.2.** *Suppose Assumption C is satisfied and $f$ is piecewise Lipschitz continuous on $[0,1]$, i.e. there exist a partition $I_1,\dots,I_k, k<\infty$, with $\bigcup_{i=1}^k I_k = [0,1]$ and $I_j \cap I_r = \emptyset$ for $j \neq r$ such that $f|_{I_j}$ is Lipschitz for all $j=1,\dots,k$. Then*

$$\int_0^1 f(x)h(x)dx = \frac{1}{n}\sum_{i=1}^n f(x_i) + O_P(n^{-1/2})\,.$$

*If additionally Assumption B is met*

$$\|\Phi f\|_2^2 = O_P(\|\Phi f\|_n^2)\,.$$

*Proof.* The proof is straightforward. For details see Boysen (2006) Lemma 7.2. □

### *4.2. Entropy results*

To show consistency of the estimates, we wish to apply results from empirical process theory. To this end, let us first introduce some additional notation (cf. van de Geer, 2000).

Given a measure $Q$, a set of $Q$-measurable functions $\mathcal{G}$ and a real number $\delta>0$, define the $\delta$-covering number $N(\delta,\mathcal{G},Q)$ as the smallest value of $N$ for which there exist functions $g_1,\dots,g_N$ such that for every $g\in\mathcal{G}$ there is a $j\in 1,\dots N$ with

$$\Big(\int (g-g_j)^2 dQ\Big)^{1/2} \le \delta\,.$$

Moreover, define the $\delta$-entropy $H$ of $\mathcal{G}$ as

$$H(\delta,\mathcal{G},Q) = \log N(\delta,\mathcal{G},Q)\,.$$

If $Q$ is the Lebesgue measure we will write $H(\delta,\mathcal{G})$ and $N(\delta,\mathcal{G})$ instead of $H(\delta,\mathcal{G},Q)$ and $N(\delta,\mathcal{G},Q)$. Given design points $x_1,\dots,x_n\in\mathbb{R}$, the empirical measure will be denoted by $Q_n = n^{-1}\sum_{i=1}^n \delta_{x_i}$. Note that $\|\cdot\|_n$ is the norm corresponding to the space $L_2(\mathbb{R},Q_n)$.

Finally, define the entropy integral

$$J(\delta,\mathcal{G},Q) := \max\Big(\delta\,,\,\int_0^\delta H^{1/2}(u,\mathcal{G},Q)du\Big)\,.$$

Note that for our purposes, the relevant quantity is the entropy of the space $\mathcal{G}_{k,R} = \{\Phi f : f \in T_{k,R}\}$. However, it is convenient to first calculate the entropy of $(T_{k,R}, \|\cdot\|_{L_2([a,b])})$ and then use Lemma 4.1 to infer on the space $\mathcal{G}_{k,R}$.

**Lemma 4.3.** *For $-\infty < a < b < \infty$ there exists a constant $C > 0$ independent of $\delta$,$k$ and $n$, such that*

$$H(\delta,(T_{k,R},\|\cdot\|_{L_2([a,b])})) \leq C(k+1)\big(1+\log\big(\frac{R(k+1)}{\delta}\big)\big)\,.$$

*Proof.* Define the sets

$$\Delta_\phi(\delta) = \Big\{-R + mc_2\delta : m = 0,\ldots,\lceil 2R(c_2\delta)^{-1}\rceil\Big\}$$

and

$$\Gamma(\delta) = \Big\{a + mc_1\delta^2 : m = 1,\ldots,\lfloor (b-a)(c_1\delta^2)^{-1}\rfloor\Big\}\,,$$

where $c_1, c_2$ will be defined later. Define the function class $\mathcal{H}(\delta)$ as

$$\begin{aligned}\mathcal{H}(\delta) \quad = \quad & \big\{g : g(x) = \sum_{i=1}^{k+1} b_i 1_{[\gamma_i-1,\gamma_i)}(x) : b_i \in \Delta_\phi(\delta), i = 1,\ldots,k+1,\\ & \gamma_0 = a, \gamma_{k+1} = b, \gamma_i \in \Gamma(\delta), \gamma_i < \gamma_{i+1}, i = 1,\ldots,k\big\}\,.\end{aligned}$$

Now for $g_0 \in T_{k,R}$ we can choose $g \in \mathcal{H}(\delta)$ such that $\mathrm{d}(\mathcal{J}(g),\mathcal{J}(g_0)) \leq c_1\delta^2/2$, and that for any $x \in [a,b]$ with $\mathrm{d}(x,\mathcal{J}(g)) > c_1\delta^2/2$ we have $(g_0(x) - g(x))^2 \leq c_2^2\delta^2/4$. Since $g_0$ has $k$ jumps between $a$ and $b$ we get

$$\|g_0 - g\|^2_{L_2([a,b])} \leq (b-a)c_2^2\,\frac{\delta^2}{4} + k(2R)^2 c_1\,\frac{\delta^2}{2}\,.$$

Choosing $c_1 = (4kR^2)^{-1}$ and $c_2 = (b-a)^{-1/2}$ gives $\|g_0 - g\|_2 \leq \delta$. Hence $\mathcal{H}(\delta)$ is an $\delta$-covering of $(T_{k,R}, \|\cdot\|_{L_2([a,b])})$. Since

$$\#\mathcal{H}(\delta) = \left\lceil \frac{2R\sqrt{b-a}}{\delta}\right\rceil^{k+1}\left\lceil\frac{(b-a)4kR^2}{\delta^2}\right\rceil^{k} = O\Big(\Big(\frac{R(k+1)}{\delta}\Big)^{3k+1}\Big)$$

the claim is proved. ∎

Lemma 4.3 directly gives that $(T_{k,R}, \|\cdot\|_{L_2([a,b])})$ is totally bounded for $-\infty < a < b < \infty$. Note that $(T_{k,R}, \|\cdot\|_{L_2([a,b])})$ also contains functions with less than $k$ jumps and hence is closed. Consequently, it is compact.

**Corollary 4.4.** *The space $(T_{k,R}, \|\cdot\|_{L_2([a,b])})$ is compact for all $a, b$ satisfying $-\infty < a < b < \infty$.*

We will now use the assumptions on the operator $\Phi$ or, to be more precise, Lemma 4.1, to deduce bounds on the entropy of the space

$$\mathcal{G}_{k,R}(\Phi) := \{\Phi g : g \in T_{k,R}\}.$$

**Corollary 4.5.** *Assume $\Phi$ satisfies Assumption B. There exists a constant $C_2$ independent of $n$,$k$ and $R$ such that*

$$H(\delta, \mathcal{G}_{k,R}(\Phi), Q_n) \leq C_2(k+1)\big(1 + \log\big(\frac{R(k+1)}{\delta}\big)\big).$$

*Proof.* By Lemma 4.1, (i) there exist $-\infty < a < b < \infty$ and $0 < C_0 < \infty$ such that

$$\|\Phi f - \Phi g\|_n \leq C_0 \|f - g\|_{L_2([a,b])}$$

for $f, g \in T_k$. Assume $\mathcal{H}(\delta)$ is a $\delta$-covering of $(T_{k,R}, \|\cdot\|_{L_2([a,b])})$ for every $\delta > 0$. Then $\mathcal{H}(\delta/C_0)$ is a $\delta$-covering of $\mathcal{G}_K(R)$. Consequently, the claim follows from Lemma 4.3. □

Again, this implies that the space $\mathcal{G}_{k,R}(\Phi)$ equipped with the empirical norm $\|\cdot\|_n$ is compact. Consequently the functional $\|\cdot - Y\|_n$ has a minimizer in $\mathcal{G}_{k,R}(\Phi)$ for every $k$. As $\lambda_n J_\#(\cdot)$ is strictly increasing in the number of jumps for every $\lambda_n > 0$ this implies the following lemma.

**Lemma 4.6.** *For each $\lambda_n > 0$ the functional $\|\cdot - Y\|_n + \lambda_n J_\#(\cdot)$ has a minimizer in $\mathcal{G}_{\infty,R}(\Phi)$.*

### 4.3. Consistency

To deduce consistency of the jump estimates from the $L_2$ consistency of the function estimator, a result on the dependency of $\mathrm{d}(\mathcal{J}(f),\mathcal{J}(g))$ on the $L_2$ distance of $f$ and $g$ is needed. This is given by the following lemma.

**Lemma 4.7.** *Assume $f,g \in T_\infty$. Then*

$$\mathrm{d}(\mathcal{J}(f),\mathcal{J}(g)) \leq \frac{4\|f-g\|_2^2}{(\min\{|f(t_+)-f(t_-)| : t \in \mathcal{J}(f)\})^2}\,.$$

*Proof.* Let $\tau \in \mathcal{J}(f)$ and $\gamma \in \mathcal{J}(g)$, such that $|\tau-\gamma| = \mathrm{d}(\mathcal{J}(f),\mathcal{J}(g))$. Then

$$\|f-g\|_2^2 \geq |\tau-\gamma|\Big(\frac{\min\{|f(t_+)-f(t_-)| : t \in \mathcal{J}(f)\}}{2}\Big)^2,$$

which proves the assertion. □

In order to show consistency of $\hat{f}_n$, we first prove the consistency of $\Phi\hat{f}_n$. To this end we require the following result which follows directly from the proof of Theorem 4.8, page 56 in van de Geer (2000).

**Lemma 4.8.** *Assume $\varepsilon_1,\ldots,\varepsilon_n$ are i.i.d. with mean zero and $\mathrm{E}(\varepsilon_1^2) = \sigma^2 < \infty$. Set $\mathcal{G}_n(R) = \{g \in \mathcal{G} : \|g\|_n \leq R\}$ and suppose that*

$$\frac{1}{n}\,H(\delta,\mathcal{G}_n(R),Q_n) \to 0 \quad \textit{for all} \quad \delta > 0, R > 0.$$

*Then*

$$\sup_{g\in\mathcal{G}_n(R)} \big|\langle\varepsilon,g\rangle_n\big| = \sup_{g\in\mathcal{G}_n(R)} \Big|\frac{1}{n}\sum_{i=1}^n \varepsilon_i g(x_i)\Big| = o_P(1)$$

*for every $R > 0$.*

Now we are able to prove consistency of $\hat{f}_n$.

**Lemma 4.9.** *Suppose the Assumptions A, B and C are met. Then $\Phi^{-1}$ is continuous as mapping from $\{\Phi f : f \in T_{k,R}\} \subset L_2([0,1])$ to the space $(T_{k,R},\|\cdot\|_{L_2([-\epsilon,1+\epsilon])})$ for any $k \in \mathbb{N}, R > 0$. Moreover $\|\Phi f - \Phi\hat{f}_n\|_2 = o_P(1)$ and consequently*

$$\|f-\hat{f}_n\|_2 = o_P(1)\,. \tag{14}$$

*Proof.* Use (9) and $Y = \Phi f + \varepsilon$ to obtain

$$\begin{aligned} \|\Phi \hat{f}_n - \Phi f\|_n &\leq 2\langle \Phi(\hat{f}_n - f), \varepsilon_n \rangle_n + o(n^{-1}) \\ &\leq 2 \sup_{g \in \mathcal{G}_{2k,2R}(\Phi)} \big|\langle g, \varepsilon_n \rangle_n\big| + o(n^{-1}), \end{aligned}$$

since $f - \hat{f}_n \in T_{2k,2R}$. By Corollary 4.5

$$n^{-1} H(\delta, \mathcal{G}_{2k,2R}(\Phi), Q_n) \to 0 \quad \text{for all} \quad \delta > 0 .$$

Hence Lemma 4.8 gives

$$\sup_{g \in \mathcal{G}_{2k,2R}(\Phi)} |\langle g, \varepsilon_n \rangle_n| = o_P(1) .$$

This proves $\|\Phi f - \Phi \hat{f}_n\|_n = o_P(1)$. Application of Lemma 4.2 yields

$$\|\Phi f - \Phi \hat{f}_n\|_2 = o_P(1) . \tag{15}$$

Note that $\Phi$ is a linear operator and $f - \hat{f} \in T_{2k,2R}$. By Corollary 4.4 the space $(T_{2k,2R}, \|\cdot\|_{L_2([-\epsilon,1+\epsilon])})$ is compact for each $\epsilon \geq 0$. Lemma 4.1, (iii) and (ii) yield that there exists an $\epsilon \geq 0$ such that the map

$$\Phi : (\ T_{2k,2R}, \|\cdot\|_{L_2([-\epsilon,1+\epsilon])}\ ) \to L_2([0,1])$$

is continuous and one-to-one.

The inverse of a continuous injective mapping $f$ restricted to the image $f(\Omega)$ is continuous if $\Omega$ is compact. This gives continuity of $\Phi^{-1}$ as mapping from $\{\Phi f : f \in T_{2k,2R}\} \subset L_2([0,1])$ to $(T_{2k,2R}, \|\cdot\|_{L_2([-\epsilon,1+\epsilon])})$. Hence, $\|\Phi f\|_2 \to 0$ implies $\|f\|_{L_2([-\epsilon,1+\epsilon])} = \|\Phi^{-1}\Phi f\|_{L_2([-\epsilon,1+\epsilon])} \to 0$ for $f \in T_{2k,2R}$. Consequently (15) implies

$$\|f - \hat{f}\|_2 \leq \|f - \hat{f}\|_{L_2([-\epsilon,1+\epsilon])} = o_P(1). \qquad \square$$

This allows us to infer the consistency of the parameter estimates. The following corollary is a direct consequence of Lemma 4.7 and 4.9.

**Corollary 4.10.** *Suppose the prerequisites of Lemma 4.9 are met. In this case*

$$\mathrm{d}(\mathcal{J}(f), \mathcal{J}(\hat{f}_n)) = o_P(1)\,,$$

*as well as* $\#\mathcal{J}(f) = \#\mathcal{J}(\hat{f}_n)$*. Moreover, if* $f$ *is given by (2) and* $\hat{f}_n$ *by (10), we have for the estimates* $\hat{b}_i$ *of the levels* $b_i$ *that*

$$\max_{i=1,\ldots,k+1} |\hat{b}_i - b_i| = o_p(1)\,.$$

### *4.4. Asymptotic normality*

To show asymptotic normality for M-estimators, it is common to assume existence of the derivative of the function which is minimized. However, as $\phi$ is allowed to have discontinuities, a less restrictive result is needed.

As discussed in Chapter 5.3 of van der Vaart (1998) it is sufficient to assume existence of a second order Taylor-type expansion. Following this idea, the next theorem gives the asymptotic normality of the minimizer of a process $Z_n(\theta)$, provided it allows for a certain expansion. It is similar to Theorem 5.23 of van der Vaart (1998), but also covers the case of non i.i.d. random variables, which is required for the fixed design.

**Theorem 4.11.** *Assume* $\Theta \subset \mathbb{R}^d$ *is open and* $\theta_0 \in \Theta$*. Let* $(Z_n(\theta))_{\theta\in\Theta}$ *be a stochastic process. Assume there exists a sequence of random variables* $(W_n)_{n\in\mathbb{N}} \subset \mathbb{R}^d$ *and a positive definite matrix* $V \in \mathbb{R}^{d\times d}$ *such that*

$$Z_n(\theta_0 + \Delta) = Z_n(\theta_0) - 2n^{-1/2} W_n^t \Delta + \Delta^t V \Delta + R_n(\Delta) \tag{16}$$

*with*

$$\sup_{\|\Delta\|\leq\delta} \frac{R_n(\Delta)}{\|\Delta\|^2 + n^{-1}} \xrightarrow{p} 0 \qquad \text{as} \quad n \to \infty\,, \delta \to 0\,, \tag{17}$$

*as well as*

$$W_n \xrightarrow{\mathcal{L}} N(0, \Gamma)\,.$$

*If $\hat{\theta}_n$ is a consistent estimator of $\theta_0$ and $\hat{\theta}_n$ is an approximate minimizer of $Z_n$, i.e.*

$$\|\hat{\theta}_n - \theta_0\| = o_P(1) \qquad \textit{and} \qquad Z_n(\hat{\theta}_n) \leq \inf_{\theta \in \Theta}(Z_n(\theta)) + o_P(n^{-1}),$$

*then*

$$\sqrt{n}(\hat{\theta}_n - \theta_0) = V^{-1}W_n + o_P(1).$$

*Proof.* The proof is straightforward and similar to the case when the second derivatives exist. For details see Boysen (2006) Theorem 7.12. ☐

**A second order expansion for the minimized process** To derive an expansion of type (16) for the problem in (9), let us first introduce some notation. For $b, \tilde{b} \in \mathbb{R}^{k+1}$ and $\tau, \tilde{\tau} \in \Gamma_k$ set

$$g(x, b, \tau) = \sum_{j=1}^{k+1} b_j \Phi 1_{[\tau_{j-1}, \tau_j)}(x).$$

and

$$Z_n(\tilde{b}, \tilde{\tau}) = \frac{1}{n}\sum_{i=1}^{n}\Big(g(x_i, b, \tau) + \varepsilon_i - g(x_i, \tilde{b}, \tilde{\tau})\Big)^2. \tag{18}$$

Assume that $f$ and the estimate $\hat{f}_n$ as defined by (9) are given by

$$f(x) = \sum_{i=1}^{k+1} b_i \Phi 1_{[\tau_{i-1}, \tau_i)}(x) \qquad \text{and} \qquad \hat{f}_n(x) = \sum_{i=1}^{k+1} \hat{b}_i \Phi 1_{[\hat{\tau}_{i-1}, \hat{\tau}_i)}(x),$$

respectively. By definition of $Z_n(\tilde{b}, \tilde{\tau})$ it is clear that

$$Z_n(\hat{b}, \hat{\tau}) \leq \min_{(\tilde{b}, \tilde{\tau}) \in [-R,R]^{k+1} \times \Gamma_k} Z_n(\tilde{b}, \tilde{\tau}) + o(n^{-1}). \tag{19}$$

To obtain an expansion for $Z_n(\tilde{b}, \tilde{\tau})$, first examine the difference $g(x, b, \tau) - g(x, \tilde{b}, \tilde{\tau})$.

**Lemma 4.12.** *Suppose Assumption B is satisfied and $\nu(x)$ is given by (12). Define $\Delta$ by*

$$\Delta = \big(\tilde{b}_1 - b_1, \tilde{\tau}_1 - \tau_1, \tilde{b}_2 - b_2, \tilde{\tau}_2 - \tau_2, \ldots, \tilde{\tau}_k - \tau_k, \tilde{b}_{k+1} - b_{k+1}\big)^t. \tag{20}$$

*Then*

$$\begin{aligned}
&g(x,b,\tau) - g(x,\tilde{b},\tilde{\tau}) \\
&= \sum_{j=1}^{k+1} b_j \Phi 1_{[\tau_{j-1},\tau_j]}(x) - \tilde{b}_j \Phi 1_{[\tilde{\tau}_{j-1},\tilde{\tau}_j]}(x) \\
&= -\Delta^t \nu(x) + O(\|\Delta\|^2) + \sum_{i=1}^{k} O(\|\tau - \tilde{\tau}\|) 1_{[x-\tau_i, x-\tilde{\tau}_i]\cap \mathcal{J}(\phi)\neq\emptyset}\,.
\end{aligned}$$

Note that $[x-\tau_i, x-\tilde{\tau}_i]\cap\mathcal{J}(\phi)\neq\emptyset$ means that $\phi$ has a discontinuity in the interval with endpoints $x-\tau_i$ and $x-\tilde{\tau}_i$.

*Proof of Lemma 4.12.* Remember $\#\mathcal{J}(\phi)<\infty$ and $\|\phi\|_\infty<\infty$.

First assume that $\tilde{\tau}_j \geq \tau_j$ and $\phi$ is continuous on $[x-\tilde{\tau}_j, x-\tau_j]$, i.e. $\mathcal{J}(\phi)\cap[x-\tilde{\tau}_j, x-\tau_j]=\emptyset$. Then for all $y\in[x-\tilde{\tau}_j, x-\tau_j]$ we have $\phi(x-y)-\phi(x-\tau_j) = O(|y-\tau_j|)$. This leads to

$$\begin{aligned}
&\Phi 1_{[\tau_{j-1},\tau_j)}(x) - \Phi 1_{[\tau_{j-1},\tilde{\tau}_j)}(x) = -\int_{\tau_j}^{\tilde{\tau}_j} \phi(x-y)dy \\
&\quad= -(\tilde{\tau}_j-\tau_j)\phi(x-\tau_j) - \int_{\tau_j}^{\tilde{\tau}_j}(\phi(x-y)-\phi(x-\tau_j))dy \\
&\quad= (\tau_j-\tilde{\tau}_j)\phi(x-\tau_j) - O(1)\int_{\tau_j}^{\tilde{\tau}_j}|y-\tau_j|dy \\
&\quad= (\tau_j-\tilde{\tau}_j)\phi(x-\tau_j) + O((\tau_j-\tilde{\tau}_j)^2)\,.
\end{aligned}$$

If $\phi$ has a discontinuity in $[x-\tilde{\tau}_j, x-\tau_j]$, then

$$\begin{aligned}
\Phi 1_{[\tau_{j-1},\tau_j)}(x) - \Phi 1_{[\tau_{j-1},\tilde{\tau}_j)}(x) &= (\tau_j-\tilde{\tau}_j)\phi(x-\tau_j) + \int_{\tau_j}^{\tilde{\tau}_j} O(\|\phi\|_\infty)dy \\
&= (\tau_j-\tilde{\tau}_j)\phi(x-\tau_j) + O(|\tau_j-\tilde{\tau}_j|)\,.
\end{aligned}$$

The same holds for $\tilde{\tau}_j<\tau_j$. Note that $1_{[x-\tau_j,x-\tilde{\tau}_j]\cap\mathcal{J}(\phi)\neq\emptyset}$ is one if and only if $\phi$ has a discontinuity in $[x-\tilde{\tau}_j, x-\tau_j]$. Consequently,

$$\begin{aligned}
\Phi 1_{[\tau_{j-1},\tau_j)}(x) - \Phi 1_{[\tau_{j-1},\tilde{\tau}_j)}(x) &= (\tau_j-\tilde{\tau}_j)\phi(x-\tau_j) + O((\tau_j-\tilde{\tau}_j)^2) + \\
&\qquad O(|\tau_j-\tilde{\tau}_j|)1_{[x-\tau_j,x-\tilde{\tau}_j]\cap\mathcal{J}(\phi)\neq\emptyset}\,.
\end{aligned}$$

Similarly,

$$\Phi 1_{[\tau_{j-1},\tau_j)}(x) - \Phi 1_{[\tilde\tau_{j-1},\tau_j)}(x) = (\tilde\tau_{j-1} - \tau_{j-1})\phi(x - \tau_{j-1}) + O((\tau_{j-1} - \tilde\tau_{j-1})^2) + O(|\tau_{j-1} - \tilde\tau_{j-1}|) 1_{[x-\tau_{j-1},x-\tilde\tau_{j-1}]\cap\mathcal{J}(\phi)\neq\emptyset}\,.$$

Remember $\tau_0 = \tilde\tau_0$ and $\tau_{k+1} = \tilde\tau_{k+1}$, combine the preceding results to obtain

$$\begin{aligned}
&\sum_{j=1}^{k+1}\Big(b_j \Phi 1_{[\tau_{j-1},\tau_j]}(x) - \tilde b_j \Phi 1_{[\tilde\tau_{j-1},\tilde\tau_j]}(x)\Big)\\
&= \sum_{j=1}^{k+1}\Big((b_j - \tilde b_j)\Phi 1_{[\tau_{j-1},\tau_j]}(x) + \tilde b_j\big(\Phi 1_{[\tau_{j-1},\tau_j]}(x) - \Phi 1_{[\tau_{j-1},\tilde\tau_j]}(x)\big) +\\
&\tilde b_j\big(\Phi 1_{[\tau_{j-1},\tilde\tau_j]}(x) - \Phi 1_{[\tilde\tau_{j-1},\tilde\tau_j]}(x)\big)\Big)\\
&= \sum_{j=1}^{k+1}\Big((b_j - \tilde b_j)\Phi 1_{[\tau_{j-1},\tau_j]}(x) + \tilde b_j(\tau_j - \tilde\tau_j)\phi(x-\tau_j) + O((\tau_j - \tilde\tau_j)^2) +\\
&O(|\tau_j - \tilde\tau_j|)1_{[x-\tau_j,x-\tilde\tau_j]\cap\mathcal{J}(\phi)\neq\emptyset} + \tilde b_j(\tilde\tau_{j-1} - \tau_{j-1})\phi(x - \tau_{j-1}) +\\
&O((\tau_{j-1} - \tilde\tau_{j-1})^2) + O(|\tau_{j-1} - \tilde\tau_{j-1}|)1_{[x-\tau_{j-1},x-\tilde\tau_{j-1}]\cap\mathcal{J}(\phi)\neq\emptyset}\Big)\,.
\end{aligned}$$

By $\tilde b_j(\tau_j - \tilde\tau_j) = b_j(\tau_j - \tilde\tau_j) + O(\|b - \tilde b\|\,\|\tau - \tilde\tau\|)$, this gives

$$\begin{aligned}
&g(x, b, \tau) - g(x, \tilde b, \tilde\tau) =\\
&\sum_{j=1}^{k+1}(b_j - \tilde b_j)\Phi 1_{[\tau_{j-1},\tau_j]}(x) + \sum_{j=1}^{k}(\tau_j - \tilde\tau_j)(b_j - b_{j+1})\phi(x - \tau_j) +\\
&O(\|\tau - \tilde\tau\|^2) + O(\|b - \tilde b\|\,\|\tau - \tilde\tau\|) + \sum_{j=1}^{k} O(\|\Delta\|)1_{[x-\tau_i,x-\tilde\tau_i]\cap\mathcal{J}(\phi)\neq\emptyset}\,.
\end{aligned}$$

Since $O(\|b - \tilde b\|\,\|\tau - \tilde\tau\|) = O(\|\Delta\|^2)$ this proves the claim. □

**Lemma 4.13.** *Suppose the Assumptions A, B and C are met. Then the process $Z_n(\tilde b, \tilde\tau)$ allows an expansion of type (16), namely*

$$Z_n(\tilde b, \tilde\tau) = Z_n(b, \tau) + 2n^{-1/2}W_n^t\Delta + \Delta^t V\Delta + R_n(\Delta)\,,$$

*where $R_n$ satisfies condition (17), $\Delta$ is given by (20) and $V$ is the $(2k+1)\times(2k+1)$ matrix defined by (13). Moreover*

$$W_n \xrightarrow{\mathcal{L}} N(0, \mathrm{E}(\varepsilon_1^2)V)\,.$$

Before we give the proof, we need the following result on the number of design points contained in a sequence of intervals.

**Lemma 4.14.** *If the design points $x_1, \dots, x_n$ satisfy Assumption C, then for any two sequences $a_n, b_n, n \in \mathbb{N}$ with $0 \le a_n < b_n \le 1$ we have*

$$n^{-1}\big(\#\{i : x_i \in [a_n, b_n]\}\big) = O_P(|b_n - a_n| + n^{-1/2}).$$

*Proof.* The proof is straightforward using that $H(x) = \int_0^x h(y)dy$ is strictly monotone, and that by Assumption C $H^{-1}(i/n - \delta_i) = x_{(i)}$ with $\max_{i=1,\dots,n} |\delta_i| = O_P(n^{-1/2})$. □

*Proof of Lemma 4.13.* Expand (18) to obtain

$$\begin{aligned} Z_n(\tilde{b}, \tilde{\tau}) = {} & \frac{2}{n}\sum_{i=1}^{n} \varepsilon_i \Big( g(x_i, b, \tau) - g(x_i, \tilde{b}, \tilde{\tau}) \Big) \\ & + \frac{1}{n}\sum_{i=1}^{n} \Big( g(x_i, b, \tau) - g(x_i, \tilde{b}, \tilde{\tau}) \Big)^2 + \|\varepsilon\|_n^2 . \end{aligned} \tag{21}$$

Note that the last term equals $Z_n(b, \tau)$. We will first estimate the second term of (21). Denote the points of discontinuity of $\phi$ by $\mathcal{J}(\phi) = \{\vartheta_1, \dots, \vartheta_{\#\mathcal{J}(\phi)}\}$ with $\vartheta_1 < \vartheta_2 < \ldots < \vartheta_{\#\mathcal{J}(\phi)}$. This means

$$[x - \tau_i, x - \tilde{\tau}_i] \cap \mathcal{J}(\phi) \neq \emptyset \quad \Leftrightarrow \quad \exists s : x \in [\vartheta_s - \tau_i, \vartheta_s - \tilde{\tau}_i] \, .$$

By Lemma 4.14,

$$\#\{i : x_i \in [\vartheta_s - \tau_j, \vartheta_s - \tilde{\tau}_j]\} = O_P(n|\tau_j - \tilde{\tau}_j| + n^{1/2}) \, .$$

This gives

$$\begin{aligned} \frac{1}{n}\sum_{i=1}^{n}\sum_{j=1}^{k}\sum_{s=1}^{\#\mathcal{J}(\phi)} 1_{[\vartheta_s - \tau_j, \vartheta_s - \tilde{\tau}_j]}(x_i) &= \frac{\#\mathcal{J}(\phi)}{n}\sum_{j=1}^{k} O_P(n|\tau_j - \tilde{\tau}_j| + n^{1/2}) \\ &= O_P(\|\Delta\| + n^{-1/2}) \, . \end{aligned}$$

The functions $\nu_j(x)$ are piecewise Lipschitz continuous by part (iv) of Lemma 4.1. With the help of Lemma 4.2 this gives

$$\begin{aligned}\frac{1}{n}\sum_{i=1}^{n}(\Delta^t\nu(x_i))^2 &= \frac{1}{n}\sum_{i=1}^{n}\sum_{j,r}^{2k+1}\Delta_j\Delta_r\nu_j(x_i)\nu_r(x_j)\\ &= \sum_{j,r}^{2k+1}\Delta_j\Delta_r\int_0^1\nu_j(x)\nu_r(x)h(x)dx + o_P(1)\\ &= \Delta^t V\Delta + o_P(1)\,.\end{aligned} \tag{22}$$

Use Lemma 4.12 and the results above to obtain

$$\begin{aligned}&\frac{1}{n}\sum_{i=1}^{n}(g(x_i,b,\tau)-g(x_i,\tilde{b},\tilde{\tau}))^2\\ &= \frac{1}{n}\sum_{i=1}^{n}\Big(\Delta^t\nu(x_i)+O(\|\Delta\|^2)+O(\|\Delta\|)\sum_{j=1}^{k}\sum_{s=1}^{\#\mathcal{J}(\phi)}1_{[\vartheta_s-\tau_j,\vartheta_s-\tilde{\tau}_j]}(x_i)\Big)^2\\ &= \frac{1}{n}\sum_{i=1}^{n}\big(\Delta^t\nu(x_i)+O(\|\Delta\|^2)\big)^2+O(\|\Delta\|^2)O_P(\|\Delta\|+n^{-1/2})\\ &= \Delta^tV\Delta+O_P(\|\Delta\|^3)+o_P(\|\Delta\|^2)\,,\end{aligned}$$

where $V$ is given by (13). The remainder terms clearly satisfy condition (17).

Next, examine the first term of (21). Set

$$W_n = n^{-1/2}\sum_{i=1}^{n}\varepsilon_i\nu(x_i)$$

to derive

$$\begin{aligned}&\frac{1}{n}\sum_{i=1}^{n}\varepsilon_i(g(x_i,b,\tau)-g(x_i,\tilde{b},\tilde{\tau}))\\ &= -\sum_{i=1}^{n}\frac{\varepsilon_i}{n}\Big(\Delta^t\nu(x_i)+O(\|\Delta\|^2)+O(\|\Delta\|)\sum_{j=1}^{k}\sum_{s=1}^{\#\mathcal{J}(\phi)}1_{[\vartheta_s-\tau_j,\vartheta_s-\tilde{\tau}_j]}(x_i)\Big)\\ &= -\frac{\Delta^tW_n}{\sqrt{n}}+\frac{O(\|\Delta\|^2)}{n}\sum_{i=1}^{n}\varepsilon_i+\frac{O(\|\Delta\|)}{n}\sum_{i=1}^{n}\sum_{j=1}^{k}\sum_{s=1}^{\#\mathcal{J}(\phi)}\varepsilon_i1_{[\vartheta_s-\tau_j,\vartheta_s-\tilde{\tau}_j]}(x_i)\,.\end{aligned}$$

The second term is clearly $o_P(\|\Delta\|^2)$.

To obtain an upper bound for the third term suppose $\vartheta_s-\tau_j<\vartheta_s-\tilde{\tau}_j$. Set

$$i_l(s,j)=\min\{i: x_{(i)}\geq\theta_s-\tau_j\}\quad\text{and}\quad i_u(s,j)=\max\{i:x_{(i)}<\theta_s-\tilde{\tau}_j\}\,.$$

Consequently,

$$\Big|\sum_{i=1}^{n}\varepsilon_i 1_{[\vartheta_s-\tau_j,\vartheta_s-\tilde\tau_j]}(x_i)\Big| = \Big|\sum_{i=i_l(s,j)}^{i_u(s,j)}\varepsilon_i\Big| .$$

By the law of the iterated logarithm for $\varepsilon_1, \varepsilon_2, \ldots$ i.i.d. with $\mathrm{E}(\varepsilon_1) = 0$ and $\mathrm{E}(\varepsilon_1^2) < \infty$ we have

$$\lim_{k_n\to\infty}\max_{j\in\{1,\ldots,k_n\}}(\mathrm{E}(\varepsilon_1^2)k_n\log\log k_n)^{-1/2}\Big|\sum_{i=1}^{j}\varepsilon_i\Big| = 1$$

almost surely. This implies for $\delta_n = i_u(s,j) - i_l(s,j)$ that

$$\max_{j=1,\ldots,\delta_n}\Big|\sum_{i=i_l(s,j)}^{i_u(s,j)}\varepsilon_i\Big| = O((\delta_n\log\log\delta_n)^{1/2})$$

holds almost surely. By Lemma 4.14,

$$\delta_n = \#\{i : \vartheta_s - \tau_j \le x_{(i)} < \vartheta_s - \tilde\tau_j\} = O_P(n|\tau_j - \tilde\tau_j| + \sqrt{n}) = O_P(n\|\Delta\| + \sqrt{n}) .$$

Consequently,

$$\sum_{i=1}^{n}\Big|\varepsilon_i 1_{[\vartheta_s-\tau_j,\vartheta_s-\tilde\tau_j]}(x_i)\Big| = O_P\Big(\sqrt{(n\|\Delta\| + n^{1/2})\log\log(n\|\Delta\| + n^{1/2})}\Big) .$$

The same can be shown for $\vartheta_j - \tau_j \ge \vartheta_j - \tilde\tau_j$. Since $\mathcal{J}(\phi)$ is a finite set and $k < \infty$, it follows that

$$\frac{O(\|\Delta\|)}{n}\sum_{i=1}^{n}\sum_{j=1}^{k}\sum_{s=1}^{\#\mathcal{J}(\phi)}\varepsilon_i 1_{[\vartheta_s-\tau_j,\vartheta_s-\tilde\tau_j]}(x_i) =$$
$$O(n^{-1}\|\Delta\|)O_P\Big(\sqrt{(n\|\Delta\| + n^{1/2})\log\log(n\|\Delta\| + n^{1/2})}\Big) . \qquad (23)$$

To verify condition (17) for this term, note that for $\|\Delta\| < n^{-1/2}$,

$$(23) = O_P\big(n^{-5/4}\sqrt{\log\log(n^{1/2})}\big) = o_P(n^{-1}) ,$$

and for $\|\Delta\| \ge n^{-1/2}$,

$$(23) = O_P\big(\|\Delta\|^{3/2}n^{-1/2}\sqrt{\log\log(n)}\big) = o_P(\|\Delta\|^2) .$$

This gives

$$\frac{1}{n}\sum_{i=1}^{n}\varepsilon_i(g(x_i, b, \tau) - g(x_i, \tilde b, \tilde\tau)) = -n^{-1/2}\Delta^t W_n + o_P(\|\Delta\|^2) + o_P(n^{-1}) .$$

Next, take a closer look at $W_n$. For any $a \in \mathbb{R}^{2k+1}$,

$$a^t W_n = \sum_{i=1}^n \varepsilon_i \Big( n^{-1/2} \sum_{j=1}^{2k+1} a_j \nu_j(x_i) \Big)$$

and by similar calculations as in (22)

$$\sum_{i=1}^n \Big( n^{-1/2} \sum_{j=1}^{2k+1} a_j \nu_j(x_i) \Big)^2 = \frac{1}{n} \sum_{i=1}^n (a^t \nu(x_i))^2 = a^t V a + o_P(1) \,.$$

By the central limit theorem and the Cramer-Wold device,

$$W_n \xrightarrow{\mathcal{L}} N(0, \sigma^2 V) \,,$$

where $\sigma^2 = \mathrm{E}(\varepsilon_1^2)$ and $V$ is given by (13). □

**Lemma 4.15.** *Given the Assumptions C and B, the matrix $V$ defined by (13) is positive definite.*

*Proof.* For any $\beta \in \mathbb{R}^{2k+1}$

$$\beta^t V \beta = \int \Big( \sum_{i=1}^{2k+1} \beta_i \nu_i(x) \Big)^2 h(x) dx \ge c_l \int_0^1 \Big( \sum_{i=1}^{2k+1} \beta_i \nu_i(x) \Big)^2 dx \,.$$

Observe that by Assumption B, the functions $\nu_1, \ldots, \nu_{2k+1}$ are linearly independent as functions in $L_2([0,1])$, since $b_i - b_{i+1} \neq 0$ for all $i = 1, \ldots, k$. Consequently, for $\beta \neq 0$ we have that

$$\int_0^1 \Big( \sum_{i=1}^{2k+1} \beta_i \nu_i(x) \Big)^2 dx > 0$$

and thus $\beta^t V \beta > 0$. □

### *4.5. Proof of Theorem 3.1*

The proof of the main theorem is now a direct consequence of the results given above. Part (v) follows directly from the proof of Lemma 4.13.

**Proof of part (i)** Corollary 4.10 implies $\|\theta - \hat\theta_n\| = o_P(1)$. By relation (19) and Lemma 4.13 the assumptions of Theorem 4.11 are satisfied. The claim follows by application of this theorem.

**Proof of part (ii)** By Lemma 4.12

$$\begin{aligned}&\int_0^1 \Big(b_i\Phi 1_{[\tau_{i-1},\tau_i)}(x) - \hat{b}_i\Phi 1_{[\hat{\tau}_{i-1},\hat{\tau}_i)}(x)\Big)^2 dx\\ &= \int_0^1 \big((\theta-\hat{\theta})^t\nu(x)\big)^2 dx + O_P(\|\theta-\hat{\theta}\|^2) = O_P(n^{-1})\,,\end{aligned}$$

since $\nu(x)$ is bounded. This proves the claim.

**Proof of part (iv) and part (iii)** Note that

$$\begin{aligned}\|f-\hat{f}_n\|_2^2 &= \sum_{i=1}^{k+1}(b_i-\hat{b}_i)^2\Big(\min(\tau_i,\hat{\tau}_i)-\max(\tau_{i-1},\hat{\tau}_{i-1})\Big)+\\ &\quad \sum_{i=1}^{k}\Big(1_{\tau_i\geq\hat{\tau}_i}(b_i-\hat{b}_{i+1})^2+1_{\tau_i<\hat{\tau}_i}(b_{i+1}-\hat{b}_i)^2\Big)|\tau_i-\hat{\tau}_i|\\ &= O_P(n^{-1})O_P(1)+O_P(1)O_P(n^{-1/2}) = O_P(n^{-1/2})\,.\end{aligned}$$

This proves part (iv). Part (iii) follows by application of Lemma 4.7. □

### 4.6. Proof of Theorem 3.2

In this section we analyze the case where the number of jumps is unknown.

In order to reconstruct the number of jumps correctly, it is helpful to use a penalty function which is strictly increasing in the number of jumps. Any penalty term, which depends on the number of jumps only, is not a pseudo-norm on $T_{\infty,R}$, since $\#\mathcal{J}(\lambda f) = \#\mathcal{J}(f)$ for $\lambda \neq 0$. Hence, the standard results from empirical process theory do not apply. However, it is possible to use similar techniques in the proofs.

The fact that $\hat{f}_{\lambda_n}$ (approximately) minimizes the penalized $L_2$ functional, implies that for any $f \in T_{\infty,R}$, we get that

$$\|\Phi\hat{f}_{\lambda_n} - Y\|_n^2 + \lambda_n J_\#(\hat{f}_{\lambda_n}) \leq \|\Phi f - Y\|_n^2 + \lambda_n J_\#(f) + o(n^{-1})\,.$$

This gives

$$\begin{aligned}&\|\Phi\hat{f}_{\lambda_n} - \Phi f\|_n^2 + 2\langle\Phi\hat{f}_{\lambda_n} - \Phi f, -\varepsilon\rangle_n + \|\varepsilon\|_n + \lambda_n J_\#(\hat{f}_{\lambda_n}) \leq\\ &\qquad\qquad \|\varepsilon\|_n + \lambda_n J_\#(f) + o(n^{-1})\,,\end{aligned}$$

which yields the basic inequality

$$\|\Phi\hat{f}_{\lambda_n} - \Phi f\|_n^2 + \lambda_n J_\#(\hat{f}_{\lambda_n}) \leq 2\langle \Phi\hat{f}_{\lambda_n} - \Phi f, \varepsilon\rangle_n + \lambda_n J_\#(f) + o(n^{-1})\,. \quad (24)$$

Hence, a bound for the term $|\langle \Phi\hat{f}_{\lambda_n} - \Phi f, \varepsilon\rangle_n|$, would allow immediate conclusions on $\|\Phi\hat{f}_{\lambda_n} - \Phi f\|_n^2$ as well as $\lambda_n J_\#(\hat{f}_{\lambda_n})$.

**Theorem 4.16.** *Suppose Assumption A is met and the error satisfies (A1). Assume* $\sup_{g\in\mathcal{G}}\|g\|_n \leq R$*. There exists a constant* $C$ *depending only on Assumption (A1), such that for all* $\delta > 0$ *satisfying*

$$\sqrt{n}\delta \geq C\Big(\int_0^R H^{1/2}(u, \mathcal{G}, Q_n)du \vee R\Big) \quad (25)$$

*we have that*

$$\mathrm{P}\Big(\sup_{g\in\mathcal{G}}\Big|\frac{1}{n}\sum_{i=1}^n \varepsilon_i g(x_i)\Big| \geq \delta\Big) \leq C\exp\Big(-\frac{n\delta^2}{C^2R^2}\Big)\,. \quad (26)$$

*Proof.* See Lemma 3.2, page 29 in van de Geer (2000). □

A bound of this type can be obtained from the following exponential inequality.

**Lemma 4.17.** *Suppose Assumptions A and B are met and the error additionally satisfies (A1).*

*There exist constants* $c_1, c_2 > 0$*, such that for all* $t \geq c_1 n^{-1/2}$ *we have*

$$\mathrm{P}\Big(\sup_{f\in T_{\infty,R}} \frac{|\langle\varepsilon, \Phi f\rangle_n|}{\|\Phi f\|_n J_\#^{1/2}(f)\big(1 + \log(J_\#(f)/\|\Phi f\|_n)_+\big)} \geq t\Big) \leq c_2\exp\Big(-\frac{nt^2}{c_2^2}\Big)\,.$$

*Proof.* Set $\mathcal{G}_{k,R}(\Phi) = \{\Phi g : g \in T_{k,R}\}$. By Corollary 4.5 there exists a constant $C > 0$ independent of $u$,$k$,$R$ and $n$ such that

$$H\big(u, \mathcal{G}_{k-1,R}(\Phi), Q_n\big) \leq Ck\big(1 + \log\big(\frac{Rk}{u}\big)\big)\,.$$

Compute

$$
\begin{aligned}
&\int_0^\delta H^{1/2}\big(u, \mathcal{G}_{k-1,R}(\Phi), Q_n\big)du \\
&\leq \sqrt{Ck}\int_0^\delta \sqrt{\log\big(\frac{\exp(1)Rk}{u}\big)}du \\
&= eRk\sqrt{Ck}\int_0^{\frac{\delta}{eRk}} \sqrt{-\log(u)}du \leq eRk\sqrt{Ck}\int_0^{\frac{\delta}{eRk}} (-\log(u))du \\
&= eRk\sqrt{Ck}\Big(\frac{\delta}{eRk}\big(1-\log\big(\frac{\delta}{eRk}\big)\big)\Big) = \delta\sqrt{Ck}(2+\log(R)+\log(k\delta^{-1})) \\
&\leq C_1\delta\sqrt{k}\big(1+\log\big(\frac{k}{\delta}\vee 1\big)\big) = C_1\delta\sqrt{k}\big(1+\log\big(\frac{k}{\delta}\big)_+\big)\,,
\end{aligned}
$$

where $C_1$ is some finite constant independent of $k$ and $\delta$. By Theorem 4.16 there exists some constant $C_2$ depending on the subgaussian error condition (A1) only, such that

$$
\sqrt{n}\rho \geq C_2\Big(\int_0^\delta H^{1/2}\big(u, \mathcal{G}_{k-1,R}(\Phi), Q_n\big)du \vee \delta\Big)
$$

implies

$$
\mathrm{P}\Big(\sup_{g\in\mathcal{G}^{(n)}_{k-1,R}(\Phi,\delta)} |\langle g,\varepsilon\rangle_n| \geq \rho\Big) \leq C_2\exp\Big(-\frac{n\rho^2}{C_2^2\delta^2}\Big)\,,
$$

where $\mathcal{G}^{(n)}_{k-1,R}(\Phi,\delta) = \{g\in\mathcal{G}_{k-1,R}(\Phi) : \|g\|_n \leq \delta\}$. Consequently, for all $t \geq C_2C_1n^{-1/2}$ we have that

$$
\begin{aligned}
&\mathrm{P}\Big(\sup_{g\in\mathcal{G}^{(n)}_{k-1,R}(\Phi,\delta)} |\langle g,\varepsilon\rangle_n| \geq t\delta\sqrt{k}\big(1+\log\big(\frac{k}{\delta}\big)_+\big)\Big) \leq \\
&\qquad\qquad C_2\exp\Big(-\frac{nt^2k\big(1+\log\big(\frac{k}{\delta}\big)_+\big)^2}{C_2^2}\Big)\,.
\end{aligned}
$$

We arrive at

$$
\begin{aligned}
&\mathrm{P}\Big(\sup_{g\in\mathcal{G}_{k-1,R}(\Phi)} \frac{|\langle\varepsilon,g\rangle_n|}{\|g\|_n\sqrt{k}\big(1+\log(k/\|g\|_n)_+\big)} \geq t\Big) \\
&\leq \sum_{s=1}^\infty \mathrm{P}\Big(\sup_{g\in\mathcal{G}_{k-1,R}(\Phi,2^{1-s}R)} |\langle\varepsilon,g\rangle_n| \geq t(2^{-s}R)\sqrt{k}\big(1+\big(\log\big(\frac{k}{2^{-s}R}\big)\big)_+\big)\Big) \\
&\leq \sum_{s=1}^\infty C_2\exp\Big(\frac{-t^2nk(1+(\log(k/R)+s\log(2))_+)}{C_2^2}\Big) \\
&\leq \sum_{s=1}^\infty C_2\exp\Big(\frac{-t^2n(1+(s\log(2)-\log(R))_+)}{C_2^2}\Big)\,.
\end{aligned}
$$

Splitting this sum at $s_R := \lceil (1+\log(R))/\log(2)\rceil$ gives

$$\begin{aligned}
&\mathrm{P}\Big(\sup_{g\in\mathcal{G}_{k-1,R}(\Phi)} \frac{|\langle\varepsilon, g\rangle_n|}{\|g\|_n\sqrt{k}\big(1+\log(k/\|g\|_n)_+\big)} \ge t\Big)\\
&\le C_2\Big\lceil\frac{1+\log(R)}{\log(2)}\Big\rceil \exp\Big(\frac{-t^2n}{C_2^2}\Big) + \sum_{s=s_R}^{\infty} C_2\exp\Big(\frac{-t^2nC_3(1+s\log(2))}{C_2^2}\Big)\\
&\le C_5\exp\Big(\frac{-t^2n}{C_2^2}\Big) + \sum_{s=1}^{\infty} C_2\exp\Big(\frac{-t^2nC_4(1+s)}{C_2^2}\Big)\\
&\le C_5\exp\Big(\frac{-t^2n}{C_2^2}\Big) + \exp\Big(\frac{-t^2nC_4}{C_2^2}\Big)\int_{s=0}^{\infty} C_2\exp\Big(\frac{-t^2nC_4 s}{C_2^2}\Big)\\
&\le C_5\exp\Big(\frac{-t^2n}{C_2^2}\Big) + \frac{C_2^3}{C_4t^2n}\exp\Big(\frac{-t^2nC_4}{C_2^2}\Big) \le C_6\exp\Big(-\frac{t^2n}{C_4^2}\Big)\,.
\end{aligned}$$

Here $C_3, C_4, C_5, C_6$ are constants depending on $C_1, C_2$ and $R$ only. The last inequality holds by $t^2n \ge C_1^2C_2^2$.

Since the constant $C_6$ does not depend on $k$, the exponential inequality also holds if we additionally take the supremum over all $k$. This proves the claim. □

The above lemma yields upper bounds for the rate of $|\langle\Phi f, \varepsilon\rangle_n|$, which are stated in the subsequent corollary.

**Corollary 4.18.** *Suppose the prerequisites of Lemma 4.17 are met. Then*

$$\sup_{f\in T_{\infty,R}} |\langle\Phi f, \varepsilon\rangle_n| = \|\Phi f\|_n\sqrt{J_\#(f)}\big(1+\log(J_\#(f)/\|\Phi f\|_n)_+\big)O_P(n^{-1/2})\,.$$

*Moreover, for each $\epsilon > 0$ we have*

$$\sup_{f\in T_{\infty,R}} |\langle\Phi f, \varepsilon\rangle_n| = \|\Phi f\|_n^{1-\epsilon}(J_\#(f))^{(1+2\epsilon)/2}\; O_P(n^{-1/2})\,.$$

*Proof.* The first equation follows directly from Lemma 4.17. To show the second equation, observe that $J_\#(f) \ge 1$ and that $\sqrt{x}(1+\log(x)) \le cx^{1/2+\epsilon}$ for $x \ge 1$, $\epsilon > 0$ and $c \ge (\epsilon^{-1} \vee 1)$ Moreover, if $c$ is large enough and $x \ge 0$ then $x(1 + \log(x^{-1})) \le cx^{1-\epsilon}$. Combine these observations to derive the second equation from the first. □

Now we are in the position to prove that with probability one the penalized estimator $\hat f_{\lambda_n}$ correctly estimates the number of jumps as $n$ tends to infinity (given a proper choice of the penalty term).

*Proof of Theorem 3.2.* Application of Corollary 4.18 to (24) gives

$$\begin{aligned}\|\Phi\hat{f}_{\lambda_n} - \Phi f\|_n^2 \le \|\Phi\hat{f}_{\lambda_n} - \Phi f\|_n^{1-\epsilon} J_\#(\hat{f}_{\lambda_n} - f)^{1/2+\epsilon} O_P(n^{-1/2}) \\ + \lambda_n(J_\#(f) - J_\#(\hat{f}_{\lambda_n})) + o(n^{-1})\,,\end{aligned} \tag{27}$$

where $\epsilon$ is given by the condition $\lambda_n n^{1/(1+\epsilon)} \to \infty$.

First, assume $J_\#(\hat{f}_{\lambda_n}) \le J_\#(f)$. Then $J_\#(\hat{f}_{\lambda_n} - f)$ is bounded and (27) implies that either

$$\|\Phi\hat{f}_{\lambda_n} - \Phi f\|_n^2 = O(\lambda_n) + o(n^{-1}) \qquad \text{or} \qquad \|\Phi\hat{f}_{\lambda_n} - \Phi f\|_n^{1+\epsilon} = O_p(n^{-1/2})\,.$$

Thus, $\|\Phi\hat{f}_{\lambda_n} - \Phi f\|_n = o_P(1)$. By Lemma 4.2, this implies $\|\Phi\hat{f}_{\lambda_n} - \Phi f\|_2 = o_P(1)$. With the help of Lemma 4.7, it follows $\mathrm{d}(\mathcal{J}(\hat{f}_{\lambda_n}), \mathcal{J}(f)) = o_P(1)$, which in turn implies $J_\#(\hat{f}_{\lambda_n}) \ge J_\#(f)$ eventually.

Now assume $J_\#(\hat{f}_{\lambda_n}) \ge J_\#(f)$. Then (27) yields

$$\|\Phi\hat{f}_{\lambda_n} - \Phi f\|_n^2 \le \|\Phi\hat{f}_{\lambda_n} - \Phi f\|_n^{1-\epsilon} J_\#(\hat{f}_{\lambda_n} - f)^{1/2+\epsilon} O_P(n^{-1/2}) + o(n^{-1})\,.$$

Assume $n_k$ is a subsequence such that $\|\Phi\hat{f}_{\lambda_{n_k}} - \Phi f\|_{n_k}^{1-\epsilon} \ge c n_k^{-1/2}$ for some $c > 0$. Dividing the last equation by $\|\Phi\hat{f}_{\lambda_{n_k}} - \Phi f\|_{n_k}^{1-\epsilon}$ gives

$$\begin{aligned}\|\Phi\hat{f}_{\lambda_{n_k}} - \Phi f\|_{n_k}^{1+\epsilon} &\le J_\#(\hat{f}_{\lambda_{n_k}} - f)^{1/2+\epsilon} O_P(n_k^{-1/2}) + o(n_k^{-1/2}) \\ &= J_\#(\hat{f}_{\lambda_{n_k}} - f)^{1/2+\epsilon} O_P(n_k^{-1/2})\,.\end{aligned}$$

This yields

$$\|\Phi\hat{f}_{\lambda_{n_k}} - \Phi f\|_{n_k}^{1-\epsilon} \le J_\#(\hat{f}_{\lambda_{n_k}} - f)^{(1+\epsilon-2\epsilon^2)/(2+2\epsilon)} O_P(n_k^{-(1-\epsilon)/(2+2\epsilon)})\,.$$

Moreover, by (27)

$$\begin{aligned}&\lambda_{n_k}(J_\#(\hat{f}_{\lambda_{n_k}}) - J_\#(f)) \le \\ &\qquad O_P(n_k^{-1/2})\|\Phi\hat{f}_{\lambda_{n_k}} - \Phi f\|_{n_k}^{1-\epsilon} J_\#(\hat{f}_{\lambda_{n_k}} - f)^{1/2+\epsilon} + o(n_k^{-1})\,.\end{aligned}$$

Combine the last two equations to obtain

$$\lambda_{n_k}(J_\#(\hat{f}_{\lambda_{n_k}}) - J_\#(f)) \le O_P(n_k^{-1/(1+\epsilon)}) J_\#(\hat{f}_{\lambda_{n_k}} - f)^{(1+\epsilon-\epsilon^2)/(1+\epsilon)}\,. \tag{28}$$

Now assume $n_k$ is a subsequence such that $\|\Phi\hat{f}_{\lambda_{n_k}} - \Phi f\|_{n_k}^{1-\epsilon} < cn_k^{-1/2}$ for some $c > 0$. Application of Corollary 4.18 to (24) and the observation that $J_\#(g) \geq 1$ for all $g$ gives

$$\begin{aligned}
&\lambda_{n_k}(J_\#(\hat{f}_{\lambda_{n_k}}) - J_\#(f)) \\
&\leq O_P(n_k^{-1/2})\|\Phi\hat{f}_{\lambda_{n_k}} - \Phi f\|_{n_k}^{1-\epsilon} J_\#(\hat{f}_{\lambda_{n_k}} - f)^{1/2+\epsilon} + o(n_k^{-1}) \\
&\leq O_P(n_k^{-1}) J_\#(\hat{f}_{\lambda_{n_k}} - f)^{1/2+\epsilon} \leq O_P(n_k^{-1/(1+\epsilon)}) J_\#(\hat{f}_{\lambda_{n_k}} - f)^{(1+\epsilon-\epsilon^2)/(1+\epsilon)} .
\end{aligned}$$

As each sequence can be decomposed into a subsequence containing only elements smaller than $cn^{-1/2}$ and a subsequence containing only elements greater or equal to $cn^{-1/2}$ for some $c > 0$, we have shown that $J_\#(\hat{f}_{\lambda_n}) \geq J_\#(f)$ implies (28).

Now we show that $J_\#(\hat{f}_{\lambda_{n_k}}) - J_\#(f) \to 0$ in probability. To this end, assume there exists some subsequence $n_k$ such that

$$J_\#(\hat{f}_{\lambda_{n_k}}) - J_\#(f) \geq c > 0 . \tag{29}$$

This implies $J_\#(f) \leq J_\#(f)c^{-1}(J_\#(\hat{f}_{\lambda_{n_k}}) - J_\#(f))$ and

$$\begin{aligned}
J_\#(\hat{f}_{\lambda_{n_k}} - f) &\leq 2(J_\#(\hat{f}_{\lambda_{n_k}}) - J_\#(f)) + 2J_\#(f) \\
&\leq (2 + 2J_\#(f)c^{-1})(J_\#(\hat{f}_{\lambda_{n_k}}) - J_\#(f)) \\
&= O(1)(J_\#(\hat{f}_{\lambda_{n_k}}) - J_\#(f)) .
\end{aligned}$$

Hence

$$\begin{aligned}
&O_P(n_k^{-1/(1+\epsilon)}) J_\#(\hat{f}_{\lambda_{n_k}} - f)^{(1+\epsilon-\epsilon^2)/(1+\epsilon)} = \\
&\qquad\qquad O_P(n_k^{-1/(1+\epsilon)})\big(J_\#(\hat{f}_{\lambda_{n_k}}) - J_\#(f)\big)^{(1+\epsilon-\epsilon^2)/(1+\epsilon)} .
\end{aligned}$$

Together with (28), the assumption $\lambda_{n_k} n_k^{1/(1+\epsilon)} \to \infty$ and (29), this gives

$$0 < c^{\epsilon^2/(1+\epsilon)} \leq \big(J_\#(\hat{f}_{\lambda_{n_k}}) - J_\#(f)\big)^{\epsilon^2/(1+\epsilon)} = O_P(\lambda_{n_k}^{-1} n_k^{-1/(1+\epsilon)}) = o_P(1) ,$$

which is a contradiction and implies $J_\#(\hat{f}_n) - J_\#(f) \to 0$ in probability. Since $J_\#(f)$ and $J_\#(\hat{f}_n)$ are integers, this yields

$$\mathrm{P}\left(J_\#(\hat{f}_n) = J_\#(f)\right) \to 1 ,$$

for $n \to \infty$. This proves the claim. □

## 5. Proof of Theorem 2.1

To give the proof of Theorem 2.1, part (i) we will define the native Hilbert space $\mathcal{N}_\phi$ of a positive definite function $\phi$ and show that the elements of its dual space $\delta_x(f) = f(x)$ and $\rho_{x,y}(f) = \int_x^y f(t)dt$ are linearly independent, if $\phi$ has certain properties. Then we will deduce that the functions $\Delta_\phi(\cdot, \tau_0, \tau_1), \dots, \Delta_\phi(x, \tau_k, \tau_{k+1})$ are linearly independent.

The assumptions $\widehat{\phi}(x) \geq 0$ and (8) imply that the Fourier transform $\widehat{\phi}$ is strictly positive. This means that $\phi$ is positive definite. (For a definition and characterization of real-valued positive definite functions, compare Chapter 6 in Wendland (2005).)

For a positive definite function $\phi$ and $\Omega \subset \mathbb{R}$ let $\mathcal{N}_\phi(\Omega)$ denote the unique Hilbert space $(\mathcal{H}, \langle\cdot,\cdot\rangle_{\mathcal{H}})$ of functions $f : \Omega \to \mathbb{R}$ satisfying $f(x) = \langle f, \phi(x-\cdot)\rangle_{\mathcal{H}}$. $\mathcal{N}_\phi(\Omega)$ is called native space for $\phi$ and given by the closure of the span of the function set $\{\phi(x-\cdot) : x \in \Omega\}$ under the inner product induced by $\langle \phi(x - \cdot), \phi(y-\cdot)\rangle = \phi(x-y)$. A short introduction to native spaces along with some basic results of the theory can be found in Schaback (1999).

Denote by

$$\mathcal{S}(\mathbb{R}) = \left\{ f \in C^\infty(\mathbb{R}, \mathbb{C}) : \lim_{|x|\to\infty} |x^n f^{(m)}(x)| = 0 \text{ for all } n, m = 0, 1, 2, \dots \right\}$$

the Schwartz space, where $C^\infty(\mathbb{R}, \mathbb{C})$ is the set of smooth functions from $\mathbb{R}$ to $\mathbb{C}$. The first result is, that the native space $\mathcal{N}_\phi(\Omega)$ contains all Schwartz functions which are compactly supported in $\Omega$.

**Lemma 5.1.** *Assume $\Omega \subset \mathbb{R}$ and $\phi$ satisfies the conditions given by Theorem 2.1, part (i). Then all real Schwartz functions with support contained in $\Omega$ are elements of the native space $N_\phi(\Omega)$, this means that*

$$\left\{ f \in \mathcal{S}(\mathbb{R}) : \operatorname{supp}(f) \subset \Omega \right\} \subset \mathcal{N}_\phi(\Omega)\,.$$

*Proof.* We first proof the claim for $\Omega = \mathbb{R}$. Assume $f \in \mathcal{S}(\mathbb{R})$. Since Fourier transformation is a bijection from $\mathcal{S}(\mathbb{R})$ to $\mathcal{S}(\mathbb{R})$ $\widehat{f}$ and $\widehat{f}^{\,2}$ are also Schwartz functions. Hence for any $n_0 \in \mathbb{N}$, we can find a constant $c_1 > 0$ such that $|\widehat{f}(x)|^2 \leq c_1(1+|x|^{n_0+2})^{-1}$. By (8) there exist $c_2 > 0$ and $n_0 \in \mathbb{N}$ such that $(\widehat{\phi}(x))^{-1} \leq c_2(1+|x|^{n_0})$. We arrive at

$$\int_{\mathbb{R}} \frac{|\widehat{f}(x)|^2}{\widehat{\phi}(x)} dx \leq c_1 c_2 \int_{\mathbb{R}} \frac{1+|x|^{n_0}}{1+|x|^{n_0+2}} dx < \infty \,.$$

By Theorem 10.12 of Wendland (2005) the function $f$ is in $\mathcal{N}_\phi(\mathbb{R})$ if and only if

$$\int_{\mathbb{R}} |\widehat{f}(x)|^2/\widehat{\phi}(x) dx < \infty \,.$$

This proves the claim for $\Omega = \mathbb{R}$.

Now assume $\Omega \subset \mathbb{R}$ is arbitrary and $f \in \mathcal{S}(\mathbb{R})$ with $\operatorname{supp} f \subset \Omega$. We have shown $f \in \mathcal{N}_\phi(\mathbb{R})$. By Theorem 10.47 in Wendland (2005) for $\Omega \subset \mathbb{R}$, $f \in \mathcal{N}_\phi(\mathbb{R})$ implies $f|_\Omega \in \mathcal{N}_\phi(\Omega)$. This proves the claim. □

Note that Lemma 5.1 implies that for any interval $(a,b) \subset \Omega$ there exists some test function $\psi \in \mathcal{N}_\phi(\Omega)$ satisfying $\operatorname{supp}(\psi) = [a,b]$. One example is

$$\psi(x) = 1_{(a,b)}(x)\, \exp((x-a)^{-1} + (b-x)^{-1}) \,.$$

This observation can be used to show that point evaluation and integral mean are linearly independent as elements of the dual space of $\mathcal{N}_\phi(\Omega)$.

**Definition 5.2.** *For $\gamma \in \mathbb{R}$ and $\gamma_1, \gamma_2 \in \mathbb{R} \cup \{-\infty, \infty\}$ with $\gamma_1 \leq \gamma_2$ define the point evaluation functional $\delta_\gamma(f) = f(\gamma)$ and the functional $\rho_{\gamma_1,\gamma_2} : \mathcal{N}_\phi(\Omega) \to \mathbb{R}$ by*

$$\rho_{\gamma_1,\gamma_2}(f) := \begin{cases} \displaystyle\int_{\gamma_1}^{\gamma_2} f(x)dx & \gamma_1 \neq \gamma_2 \,, \\ f(\gamma_1) & \gamma_1 = \gamma_2 \,. \end{cases}$$

**Lemma 5.3.** *Suppose $\phi$ satisfies the conditions given by Theorem 2.1, part (i). Assume $\tau_0 < \ldots < \tau_{k+1}$, $\gamma_1 < \ldots < \gamma_r$ and there exist an $\epsilon > 0$ such that $(\tau_1 - \epsilon, \tau_k + \epsilon) \subset \Omega$ as well as $(\gamma_1 - \epsilon, \gamma_r + \epsilon) \subset \Omega$. Then the functionals* $\rho_{\tau_0,\tau_1}, \rho_{\tau_1,\tau_2}, \ldots, \rho_{\tau_k,\tau_{k+1}}, \delta_{\gamma_1}, \ldots, \delta_{\gamma_r}$ *are linearly independent as elements of the dual space* $\mathcal{N}_\phi(\Omega)'$.

*Proof.* Assume

$$\sum_{i=1}^{k+1} \alpha_i \rho_{\tau_{i-1},\tau_i}(f) + \sum_{j=1}^{r} \beta_j \delta_{\gamma_j}(f) = 0$$

for all $f \in \mathcal{N}_\phi(\Omega)$. For each $i = 1, \ldots, k+1$ we can find an interval $J_i \subset [\tau_{i-1}, \tau_i] \cap \Omega$ such that $J_i \cap \gamma_j = \emptyset$ for all $j = 1, \ldots, r$. By Lemma 5.1 we can find a test function $f_i \in \mathcal{N}_\phi(\Omega)$ with $\mathrm{supp}(f_i) \subset J_i$ and $\int_{\mathbb{R}} f_i(x)dx = 1$ for all $i = 1, \ldots, k+1$. We then have that $\rho_{\tau_{l-1},\tau_l}(f_i) = 1_{i=l}$ and $\delta_{\gamma_j}(f_i) = 0$ for all $i = 1, \ldots, k+1$ and $j = 1, \ldots, r$. This leads to

$$0 = \sum_{l=1}^{k+1} \alpha_l \rho_{\tau_{l-1},\tau_l}(f_i) + \sum_{j=1}^{r} \beta_j \delta_{\gamma_j}(f_i) = \alpha_i$$

for all $i = 1, \ldots, k+1$. Similarly we can find test functions $f_j \in \mathcal{N}_\phi(\Omega)$ with $\delta_{\gamma_j}(f_i) = 1_{i=j}$ and deduce that $\beta_j = 0$ for all $j = 1, \ldots, r$. This proves the claim. □

Finally, we can prove Theorem 2.1, part (i).

*Proof of Theorem 2.1, part (i).* Assume

$$\Big\| \sum_{i=1}^{k+1} \alpha_i \Delta_\phi(\cdot, \tau_{i-1}, \tau_i) \Big\|_2 = 0\,. \tag{30}$$

By continuity of $\phi$, $\Delta_\phi(x, \tau_{i-1}, \tau_i)$ and hence $\sum_{i=1}^{k+1} \alpha_i \Delta_\phi(x, \tau_{i-1}, \tau_i)$ are continuous functions of $x$. Consequently, (30) implies

$$0 = \sum_{i=1}^{k+1} \alpha_i \Delta_\phi(x, \tau_{i-1}, \tau_i)\,,$$

for all $x \in [0,1]$. By definition of $\Delta_\phi$ (see (7))

$$0 = \sum_{i=1}^{k+1} \alpha_i \Delta_\phi(x, \tau_{i-1}, \tau_i) = \sum_{i=1}^{k+1} \alpha_i \rho_{\tau_{i-1},\tau_i}\big(\phi(x - \cdot)\big)\,,$$

for all $x \in [0,1]$. Set $\Omega = [0,1]$. By Theorem 8 in Schaback (1999), the native space $\mathcal{N}_\phi(\Omega)$ is the closure of the span of the set of functions $\{\phi(x - \cdot) : x \in \Omega\}$. It follows that

$$0 = \sum_{i=1}^{k+1} \alpha_i \rho_{\tau_{i-1},\tau_i}(f)$$

for all $f \in \mathcal{N}_\phi(\Omega)$. By Lemma 5.3 we know that $\rho_{\tau_0,\tau_1}, \dots, \rho_{\tau_k,\tau_{k+1}}$ are linearly independent as elements of the dual space $\mathcal{N}_\phi(\Omega)'$. Consequently, $\alpha_i = 0$ for all $i = 1, \dots, k+1$, which proves the claim. □

## 6. A lower bound for estimating the jump locations

In this section we show that the obtained rate $\mathrm{d}(\mathcal{J}(\hat{f}_n), \mathcal{J}(f)) = O_P(n^{-1/2})$ is optimal in a minimax sense. To do so, we construct functions $f_0, f_{1,n}, f_{2,n}$ with $\mathrm{d}(\mathcal{J}(f_0), \mathcal{J}(f_{i,n})) = cn^{-1/2}$ for $i = 1, 2$ and some $c > 0$ to be chosen later. Given the observations

$$Y_i = g(x_i) + \varepsilon_i \quad i = 1, \dots, n$$

for $g \in \{\Phi f_0, \Phi f_{1,n}, \Phi f_{2,n}\}$ and $\varepsilon_1, \dots, \varepsilon_n$ independent and identically distributed according to $N(0, \sigma^2)$ with $\sigma^2 > 0$, we show that for any estimator, the probability to choose the true function is strictly smaller than one. Obviously it is sufficient to consider the case of a single jump with a fixed jump height.

**Lemma 6.1.** *Suppose Assumption B is met, $x_1, \dots, x_n \in [0,1]$ are arbitrary fixed design points. Moreover, assume that $\varepsilon_1, \dots, \varepsilon_n$ are independent and identically distributed according to $N(0, \sigma^2)$ with $\sigma^2 > 0$. Set $g_\tau = \Phi 1_{[\tau,\infty)}$ for $\tau \in \mathbb{R}$. Given observations*

$$Y_i = g_\tau(x_i) + \varepsilon_i \quad i = 1, \dots, n$$

*denote the corresponding probability measure by $P_\tau$. There exist constants $c, c_1 > 0$ such that*

$$\inf_{\hat{\tau}} \sup_{\tau \in \mathbb{R}} P_\tau(|\tau - \hat{\tau}| \geq cn^{-1/2}) \geq c_1 > 0\,.$$

For the proof we need the following theorem.

**Theorem 6.2.** *Suppose $M \geq 2$ and $\Theta$ contains elements $\theta_0, \theta_1, \ldots, \theta_M$ with*

$$d(\theta_j, \theta_k) \geq 2s > 0, \qquad \forall\, 0 \leq j < k \leq M\,,$$

*$P_{\theta_j} \ll P_{\theta_0}$ for all $j = 1, \ldots, M$ and*

$$\frac{1}{M}\sum_{j=1}^{M} d_K(P_{\theta_j}, P_{\theta_0}) \leq \alpha \log M\,,$$

*with $0 < \alpha < 1/10$, where $d_K$ denotes the Kullback-Leibler distance. Then*

$$\inf_{\hat\theta} \sup_{\theta\in\Theta} P_\theta(d(\hat\theta, \theta) \geq s) \geq \frac{\sqrt{M}}{1+\sqrt{M}}\Big(1 - 2\alpha - 2\sqrt{\frac{\alpha}{\log M}}\Big) > 0\,.$$

*Proof.* See Theorem 2.5, page 85 in Tsybakov (2004). □

*Proof of Lemma 6.1.* Recall that for product measures $P = \otimes_{i=1}^n P_i$ and $Q = \otimes_{i=1}^n Q_i$ we have

$$d_K(P, Q) = \sum_{i=1}^{n} d_K(P_i, Q_i)\,.$$

Note that $Y_i \sim N(g_\tau(x_i), \sigma^2)$ and denote the corresponding measures with $P_\tau^i$. By independence of the $\varepsilon_i$ the joint measure $P_\tau$ of $Y_1, \ldots, Y_n$ is given by $P_\tau = \otimes_{i=1}^n P_\tau^i$. This yields

$$d_K(P_{\tau_1}, P_{\tau_2}) = (2\sigma^2)^{-1}\sum_{i=1}^{n}(g_{\tau_1}(x_i) - g_{\tau_2}(x_i))^2\,.$$

By Assumption B the integral kernel $\phi$ is bounded in sup-norm. Calculate

$$(g_{\tau_1}(x_i) - g_{\tau_2}(x_i))^2 = \Big(\int_{[\tau_1,\tau_2]} \phi(x_i - y)dy\Big)^2 \leq (\tau_1 - \tau_2)^2\|\phi\|_\infty^2\,.$$

Consequently,

$$d_K(P_{\tau_1}, P_{\tau_2}) \leq (2\sigma^2)^{-1} n(\tau_1 - \tau_2)^2\|\phi\|_\infty^2\,.$$

Now choose some $0 < \alpha < 1/10$, set $c = (2\alpha\sigma^2/\|\phi\|_\infty^2)^{1/2}$ and choose

$$\tau_0 \in (\tau_{low} + cn^{-1/2}, \tau_{up} - cn^{-1/2})\,.$$

Set

$$\tau_1 = \tau_0 + cn^{-1/2} \qquad \text{and} \qquad \tau_2 = \tau_0 - cn^{-1/2}\,.$$

This gives

$$\frac{1}{2}\sum_{j=1}^{2} d_K(P_{\tau_j}, P_{\tau_0}) \leq (4\sigma^2)^{-1} n \sum_{j=1}^{2} (\tau_0 - \tau_j)^2 \|\phi\|_\infty^2 = \alpha \,.$$

Consequently, the assumptions of Theorem 6.2 are satisfied for $s = n^{-1/2}c/2$ and $d(\tau, \tau') = |\tau - \tau'|$. Application of this theorem gives

$$\inf_{\hat{\tau}} \sup_{\tau \in \mathbb{R}} P_\tau(|\tau - \hat{\tau}| \geq 2^{-1} c n^{-1/2}) \geq \frac{\sqrt{2}}{1+\sqrt{2}} \Big( 1 - 2\alpha - 2\sqrt{\frac{\alpha}{\log 2}} \Big) > 0 \,.$$

This proves the claim. □

Note that in the proof we used the absolute integrability and the boundedness in supremum norm of the integral kernel only.

*Proof of Theorem 3.3.* Lemma 6.1 directly implies that the jump estimator attains the minimax rate. If $f$ is a step function with known jump locations and unknown level heights $b_i$, the inverse regression model (1) reduces to a standard linear regression model. It is well known that in this setting the levels $b_i$ cannot be estimated at a rate faster than $O_P(n^{-1/2})$. Consequently, this also holds for the case of unknown jump locations. This proves Theorem 3.3. □

**Acknowledgments** The authors wish to thank L. Dümbgen, T. Hohage, R. Schaback and A. Tsybakov for interesting comments and bibliographic information.

L. Boysen gratefully acknowledges support by Georg Lichtenberg program "Applied Statistics & Empirical Methods" and DFG graduate program 1023 "Identification in Mathematical Models"; A. Munk was supported by DFG grant "Statistical Inverse Problems under Qualitative Shape Constraints" and DFG grant FOR916.

## References

Birgé, L. and Massart, P. (2006). Minimal penalties for gaussian model selection. *Probab. Theory Related Fields, to appear* .

BISSANTZ, N., DÜMBGEN, L., HOLZMANN, H. and MUNK, A. (2007). Nonparametric confidence bands in deconvolution density estimation. *J. Royal Statist. Society Ser. B.* **69** 483–506.

BOYSEN, L. (2006). *Jump estimation for noisy blurred step function.* Ph.D. thesis, Georg-August-Universiät Göttingen.
URL http://webdoc.sub.gwdg.de/diss/2006/boysen/boysen.pdf

BUTUCEA, C. and TSYBAKOV, A. (2007). Sharp optimality for density deconvolution with dominating bias. i,ii. *Theory of Probability and Its Applications, to appear* .

CAVALIER, L. and TSYBAKOV, A. (2002). Sharp adaption for inverse problems with random noise. *Prob. Theory Rel. Fields* **123** 323–354.

DÜMBGEN, L. and JOHNS, R. B. (2004). Confidence bands for isotonic median curves using sign tests. *J. Comput. Graph. Statist.* **13** 519–533.

FEDER, P. I. (1975). On asymptotic distribution theory in segmented regression problems–identified case. *Ann. Statist.* **3** 49–83.

GOLDENSHLUGER, A., JUDITSKY, A., TSYBAKOV, A. and ZEEVI, A. (2006a). Change-point estimation from indirect observations. 1. minimax complexity. Preprint.

GOLDENSHLUGER, A., JUDITSKY, A., TSYBAKOV, A. and ZEEVI, A. (2006b). Change-point estimation from indirect observations. 2. adaptiation. Preprint.

GOLDENSHLUGER, A., TSYBAKOV, A. and ZEEVI, A. (2006c). Optimal change-point estimation from indirect observations. *Ann. Statist.* **34** 350–372.

HINKLEY, D. V. (1969). Inference about the intersection in two-phase regression. *Biometrika* **56** 495–504.

KARLIN, S. and STUDDEN, W. J. (1966). *Tchebycheff systems: With applications in analysis and statistics.* Pure and Applied Mathematics, Vol. XV, Interscience Publishers John Wiley & Sons, New York-London-Sydney.

KOUL, H. L., QIAN, L. and SURGAILIS, D. (2003). Asymptotics of $M$-estimators in two-phase linear regression models. *Stochastic Process. Appl.*

**103** 123–154.

Lan, Y., Banerjee, M. and Michailidis, G. (2007). Change-point estimation under adaptive sampling. *Technical Report, Univ. of Michigan* .

Müller, H.-G. (1992). Change-points in nonparametric regression analysis. *Ann. Statist.* **20** 737–761.

Müller, H.-G. and Stadtmüller, U. (1999). Discontinuous versus smooth regression. *Ann. Statist.* **27** 299–337.

Neumann, M. H. (1997). Optimal change-point estimation in inverse problems. *Scand. J. Statist.* **24** 503–521.

Quandt, R. E. (1958). The estimation of the parameters of a linear regression system obeying two separate regimes. *J. Amer. Statist. Assoc.* **53** 873–880.

Quandt, R. E. (1960). Tests of the hypothesis that a linear regression system obeys two separate regimes. *J. Amer. Statist. Assoc.* **55** 324–330.

Raimondo, M. (1998). Minimax estimation of sharp change points. *Ann. Statist.* **26** 1379–1397.

Römer, W., Lam, Y. H., Fischer, D., Watts, A., Fischer, W. B., Göring, P., Wehrspohn, R. B., Gösele, U. and Steinem, C. (2004). Channel activity of a viral transmembrane peptide in micro-blms. *J. Am. Chem. Soc.* **49** 16267–16274.

Roths, T., Maier, D., Friedrich, C., Marth, M. and Honerkamp, J. (2000). Determination of the relaxation time spectrum from dynamic moduli using an edge preserving regularization method. *Rheol. Acta* **39** 163–173.

Sacks, J. and Ylvisaker, D. (1970). Designs for regression problems with correlated errors. III. *Ann. Math. Statist.* **41** 2057–2074.

Schaback, R. (1999). Native Hilbert spaces for radial basis functions. I. In *New developments in approximation theory (Dortmund, 1998)*, vol. 132 of *Internat. Ser. Numer. Math.* Birkhäuser, Basel, 255–282.

Schmitt, E. K., Vrouenraets, M. and Steinem, C. (2006). Channel activity of ompf monitored in nano-blms. *Biophys. J.* **91** 2163–2171.

Sprent, P. (1961). Some hypotheses concerning two phase regression lines. *Biometrics* **17** 634–645.

Tsybakov, A. B. (2004). *Introduction à l'estimation non-paramétrique*, vol. 41 of *Mathématiques & Applications (Berlin) [Mathematics & Applications]*. Springer-Verlag, Berlin.

van de Geer, S. A. (1988). *Regression analysis and empirical processes*, vol. 45 of *CWI Tract*. Stichting Mathematisch Centrum Centrum voor Wiskunde en Informatica, Amsterdam.

van de Geer, S. A. (2000). *Applications of empirical process theory*, vol. 6 of *Cambridge Series in Statistical and Probabilistic Mathematics*. Cambridge University Press, Cambridge.

van der Vaart, A. W. (1998). *Asymptotic statistics*, vol. 3 of *Cambridge Series in Statistical and Probabilistic Mathematics*. Cambridge University Press, Cambridge.

Wendland, H. (2005). *Scattered data approximation*, vol. 17 of *Cambridge Monographs on Applied and Computational Mathematics*. Cambridge University Press, Cambridge.

Yakir, B., Krieger, A. M. and Pollak, M. (1999). Detecting a change in regression: first-order optimality. *Ann. Statist.* **27** 1896–1913.

Yao, Y.-C. and Au, S. T. (1989). Least-squares estimation of a step function. *Sankhyā Ser. A* **51** 370–381.